\documentclass[11pt, letter]{article}
\usepackage[utf8]{inputenc}
\usepackage[authoryear,round]{natbib}
\usepackage{algpseudocode}
\usepackage{pdflscape}

\usepackage{algorithm}
\usepackage{amsfonts}
\usepackage{amsmath}
\usepackage{amssymb}
\usepackage{amsthm}
\usepackage{bbm}
\usepackage{geometry}
\usepackage{graphicx}
\usepackage{mathrsfs}
\usepackage{multicol} 
\usepackage{multirow}
\usepackage{paralist}
\usepackage{pgfplots}
\usepackage{setspace}
\usepackage{soul}
\usepackage{subcaption}
\usepackage{tikz}
\usepackage{verbatim}

\allowdisplaybreaks[4]

\definecolor{Red}{rgb}{0,0,0}
\newcommand{\Red}{\color{Red}}
\definecolor{DR}{rgb}{0,0,0}
\newcommand{\DR}{\color{DR}}
\definecolor{Blue}{rgb}{0,0,0}
\newcommand{\Blue}{\color{Blue}}
\definecolor{Green}{rgb}{0,0,0}
\newcommand{\Green}{\color{Green}}
\definecolor{Yellow}{rgb}{0,0,0}

\definecolor{PaleGrey}{rgb}{0,0,0}

\definecolor{safetyorange}{rgb}{0,0,0}

\def\vec#1{\mathchoice{\mbox{\boldmath$\displaystyle\bf#1$}}
{\mbox{\boldmath$\textstyle\bf#1$}}
{\mbox{\boldmath$\scriptstyle\bf#1$}}
{\mbox{\boldmath$\scriptscriptstyle\bf#1$}}}

\newcommand{\Px}{\mathbb{P}}
\newcommand{\Ex}{\mathbb{E}}
\newcommand{\Rx}{\mathbb{R}}

\newtheorem{theorem}{Theorem}[section]

\newtheorem{proposition}[theorem]{Proposition}
\newtheorem{remark}[theorem]{Remark}
\newtheorem{lemma}[theorem]{Lemma}
\newtheorem{corollary}[theorem]{Corollary}

\newtheorem{assumption}{Assumption}

\title{\textbf{\textsc{Optimal Iterative Threshold-Kernel Estimation\\of Jump Diffusion Processes}}}

\author{
Jos\'{e} E. Figueroa-L\'{o}pez
\thanks{Department of Mathematics and Statistics, Washington University in St. Louis, St. Louis, MO, 63130, USA.
{\tt figueroa-lopez@wustl.edu}.}
\and
Cheng Li
\thanks{Citadel Securities, New York, NY, 10022, USA.
{\tt Cheng.Li@citadelsecurities.com}.}
\and
Jeffrey Nisen
\thanks{Quantitative Analytics, Barclays, New York, NY, 10019, USA.
{\tt  jeffrey.nisen@barclayscapital.com}.}
}

\begin{document}

\maketitle

\abstract{
In this paper, we {\DR propose a  new threshold-kernel jump-detection} method for jump-diffusion processes, which iteratively applies thresholding and kernel methods in an approximately optimal way to achieve improved finite-sample performance. 
As in \cite{FLN2013}, we use the expected number of jump misclassifications as the objective function to optimally select the threshold parameter of the jump detection scheme. We prove that the objective function is quasi-convex and obtain a new second-order infill approximation of the optimal threshold {\DR in closed form}. The approximate optimal threshold depends not only on the spot volatility $\sigma_t$, but also 
the jump intensity and the value of the jump density at the origin. 
{Estimation} methods for these quantities are then developed, where the spot volatility is estimated by a kernel estimator with thresholding and the value of the jump density at the origin is estimated by a density kernel estimator applied to those increments deemed to {contain} jumps by the chosen thresholding criterion. Due to the interdependency between the model parameters and the approximate optimal estimators built to estimate them, a type of iterative fixed-point {\Red estimation} algorithm is developed to implement them. Simulation studies {for a prototypical stochastic volatility model,} show that it is not only feasible to implement the higher-order local optimal threshold scheme but also that this is superior to those based only on the first order approximation and/or on average values of the parameters over the estimation time period.}

\section{Introduction}

In this work, we study 
{a} jump diffusion process of the form
\[
X_{t} 
:= 
\int_{0}^{t} \gamma_u du + \int_{0}^{t} \sigma_u dW_{u} + \sum_{j=1}^{N_{t}} \zeta_{j},
\]
where {$W$ is a {Wiener} process, $N$ is an independent Poisson process with local intensity $\{\lambda_{t}\}_{t\geq{}0}$, and $\{\zeta_{j}\}_{j\geq{}1}$ are i.i.d. variables independent {of} $W$ and $N$.
With the presence of jumps}, several statistical inference problems, including volatility estimation and jump detection, can be addressed by the {thresholding approach {developed} by Mancini (2001, 2004, 2009)}. The basic idea is to {introduce a threshold tuning parameter $B$ so that} whenever the absolute value of {an} increment $\Delta X := X_{t_i} - X_{t_{i - 1}}$ exceeds $B$, we conclude that an unusual event (aka a ``jump") has happened during {\Red the interval} $(t_{i-1},t_{i}]$, based on which we can {then} proceed to estimate the volatility and other parameters. Many works have been conducted to further {extend} the threshold method to various statistical inference problems. For {an} It\^{o} semimartingale with finite or infinite jump activity, jump detection and integrated volatility estimation was studied by \cite{Mancini:2009} and \cite{Jacod2007,Jac08MPVperSM}. 
We also refer to \cite{Corsietal2010},  \cite{ait2009testing,ait2009estimating,ait2010brownian}, \cite{CLTIA}, \cite{figueroa2012statistical}, \cite{jing2012jump}, and others for further applications of the threshold method.

One of the key issues that we have to address in order to have {a} good performance of the {jump detection procedure} is {the selection of the threshold $B$}. {Ideally, we hope to select the best possible threshold {under a suitable criterion}. Such a problem was studied by Figueroa-L\'{o}pez and Nisen (2013) using  the expected number of jump {\Red misclassifications} {as the estimation loss function and, more recently, by \cite{FLM2017} using the mean-square error of the threshold realized quadratic variation}. Under the assumption of zero drift, constant volatility {$\sigma$}, and constant jump intensity, \cite{FLN2013} showed that} the {first-order} approximation of the optimal threshold {is given by} $\sqrt{3\,{\sigma^{2}}{h} \log(1/ {h})}$ (cf. {Theorems 4.2 and 4.3 therein), when ${h}$, the time {span} between observations, shrinks to $0$ (i.e., infill or high-frequency asymptotics)}. {\Blue Based on this result,  \cite{FLN2013} proposed a method to estimate time-dependent deterministic volatilities and, by simulation, showed that its performance is good for smooth volatilities.} 
In this work, we generalize this framework {in three directions}. 
{We first prove that the loss function is {quasi-convex} and admits a global minimum in the more general case of non-homogeneous drift, volatility, and jump intensity}. {\Blue A simpler version of this result was stated without proof in \cite{FLN2013}}. We then proceed to obtain  a second-order asymptotic approximation of the optimal localized threshold, {\DR in closed form,} which {depends on the spot volatility $\sigma_t$, the local jump intensity $\lambda_{t}$, and the value of the jump density at the origin. We find out that, as expected, 
if} the spot volatility is high, then it is more preferable to have a larger threshold. {However, when} the jump intensity or the jump density at the origin is large,  the possibility of having {smaller} jumps is higher, which favors a smaller threshold {to {detect} such jumps}. Although {an explicit formula for the second-order approximation is derived}, the method is not feasible unless we are able to estimate all the unknown parameters {appearing in this formula}: the spot volatility, the jump intensity, and the jump density at the origin. To this end, we {apply} kernel estimation techniques, as described below, to devise feasible plug-in type estimators for the optimal threshold.

{Kernel estimation} has a long history and has been applied to a large range of statistical problems.  {In our work, we use it} to estimate the jump density at the origin. The problem we are facing differs from the usual {density kernel estimation} in several ways. Firstly, the data we have {is} contaminated by {noise}, and to make things even worse, part of the data may not contain {any information at all about} the density we {want} to estimate. Moreover, due to the usage of {a} threshold, the data we have {is at best {drawn} from a truncated distribution and, the point at which we hope to estimate the density, is not even} inside the support of the truncated data. Due to these reasons, we have to adjust the {standard method of kernel density estimation} and select the threshold appropriately so that we can get a satisfactory estimation of the jump density {at the origin}. It turns out that the optimal threshold that we should use in such a situation is larger than the one we use for optimal jump detection {(see} Section \ref{sec:jump_density_estimation} for the intuition behind this).

Another quantity we have to estimate is the spot volatility, which can also be estimated by the kernel estimator. One earlier research on this topic is \cite{foster1994continuous}, where {a} rolling window estimator is analyzed, which is similar to the idea of the kernel estimation {with a uniform kernel}. The kernel-based estimation of the spot volatility, \emph{with general kernel}, was studied by \cite{fan2008spot}, \cite{kristensen2010nonparametric}, {\cite{mancini2015estimation} and,} more recently,  Figueroa-L\'{o}pez and Li (2017). See also the excellent monographs  of \cite[Ch. 13]{JacodProtter} and \cite[Ch. 8]{JacodAitSahaliaBook} for a {general treatment of the problem of spot volatility estimation of It\^o semimartigales via} uniform kernels (though Remark 8.10 in \cite{JacodAitSahaliaBook} also briefly mentions the case of a general kernel with support on $[0,1]$). One of the key issues related to kernel estimators of spot volatility is how to select the bandwidth. \cite{kristensen2010nonparametric} proposed a leave-one-out cross-validation method, which is a general method, but suffers from the loss of accuracy and computational inefficiency. In {this} work, we {adapt and extend the approach of Figueroa-L\'{o}pez and Li (2017) by applying a threshold-kernel estimator of the spot volatility rather than just kernel estimation.}  The leading order {terms} of the MSE of the estimator {are explicitly derived}, based on which {we} {propose a procedure for optimal} bandwidth and kernel selection. The {CLT of the estimation error is also given.}

{As explained above,} the approximated optimal threshold depends on the spot volatility, jump intensity, and the {value of the jump density at the origin}, while {the approximated optimal} estimators of these three quantities depend on the threshold. Such {an interdependency} immediately suggests an iterative algorithm that starts with an initial guess of these parameters and gradually {converges} to a fixed point result. Due to the {nature} of the threshold estimator, the result is purely determined by whether the absolute value of each data increment exceeds the threshold, so we can conclude convergence without any ambiguity based on whether each data increment is included by the threshold or not.

The rest of the paper is organized as follows. {Section 
\ref{sec:optimal_threshold_assumption} introduces the framework and {assumptions}. In Section \ref{section:optimal_thresholding}}, we analyze the optimal threshold and obtain the second order approximation {thereof}. {The {bias} and variance of the estimator are derived in Section 
\ref{BiasVarianceSect}}. In Section \ref{sec:jump_density_estimation}, we consider the kernel estimation of the jump density at the origin. The threshold-kernel estimation of the spot volatility is studied in Section \ref{sec:threshold_kernel_spot_vol}. The three estimators are {then} combined into an iterative algorithm {presented} in Section \ref{sec:implementation}. Finally, {the performance of the proposed methods are analyzed through several {simulations} in Section \ref{sec:monte_carlo_study}. {Conclusions and some thoughts about} future work are provided in Section \ref{sec:conclusion}. The proofs of the main results are deferred to an Appendix section.}

\section{The Optimal Threshold {of TRV}}\label{sec:optimal_threshold}
In this section we extend the modelling framework and optimal thresholding results {of \cite{FLN2013}}. Specifically, we will allow non-constant drift, volatility, and {\Red jump} intensity, though we keep the jump density {constant through time}. In the first subsection, we {introduce} all the assumptions that we need for the optimal threshold results. However, we temporarily set {the} drift, volatility, and intensity to be deterministic, {which would subsequently be relaxed} when we discuss the kernel threshold estimation of spot volatility. All the results can be generalized to {\Blue stochastic drift and volatility, and doubly stochastic Poisson process $N$, as long as we assume that the Brownian motion and jumps of the semimartingale are independent from all these processes, since we can always condition on the paths of the drift, the volatility, and the jump intensity of $N$}. {\DR It is also important to point out that,} {though our results in this section are derived under the just mentioned {\Blue independence} assumption, our simulation experiments show that this is not essential as the proposed estimators perform well under prototypical stochastic volatility {\Blue models} with leverage.}
\subsection{The Framework and Assumptions}\label{sec:optimal_threshold_assumption} 
{Throughout}, 
{we consider an 
It\^{o} semimartingale of the form:}
\begin{equation}
X_{t} := \left(\int_{0}^{t} \gamma_u du + \int_{0}^{t} \sigma_u dW_{u} \right)+ \sum_{j=1}^{N_{t}} \zeta_{j} =: X^{c}_{t} + J_{t}\label{eq:jump_diffussion_model},
\end{equation}
where {$W=\{W_t\}_{t\geq{}0}$ is a Wiener process,} $\{\zeta_{j}\}_{j\geq{}1}$ are i.i.d.~variables with density $f$, ${N=\{ N_{t}\}_{t \geq 0}}$ is a non-homogeneous Poisson process with {intensity} function $\{\lambda_t\}_{t \geq 0}$, and the continuous component $\{X^c_t\}_{t \geq 0}$ and jump component $\{J_t\}_{t \geq 0}$ are independent. {\Blue The processes $\gamma$ and $\sigma$ satisfy standard conditions for the integrals in (\ref{eq:jump_diffussion_model}) to be well-defined. In this section, we shall additionally assume the following conditions on $\gamma$, $\sigma$, and $\lambda$\footnote{{\Blue In Section \ref{sec:threshold_kernel_spot_vol}, we will consider stochastic processes $\gamma$ and $\sigma$.}}}:
\begin{assumption}\label{assumption:boundedness_of_intensity}
The functions $\gamma:[0,\infty)\to\mathbb{R}$, $\sigma:[0,\infty)\to\mathbb{R}^{+}$, and $\lambda:[0,\infty)\to\mathbb{R}^{+}$ are deterministic such that, for any given fixed {\Blue $t > 0$,
\begin{equation}\label{BdnessOnsigmaetc}
\begin{split}
	&\underline{\sigma}_{t}:=\displaystyle{\inf_{0 \leq s \leq t}}  \sigma_s > 0, \quad\bar{\sigma}_{t}:=\displaystyle{\sup_{0 \leq s \leq t}}  \sigma_s <\infty,\\
	&\underline{\gamma}_{t}:=\displaystyle{\inf_{0 \leq s \leq t}}  \gamma_{s} > 0, \quad\bar{\gamma}_{t}:=\displaystyle{\sup_{0 \leq s \leq t}}  \gamma_{s} <\infty,\\
	&\underline{\lambda}_{t}:=\displaystyle{\inf_{0 \leq s \leq t}}  \lambda_{s}> 0, \quad\bar{\lambda}_{t}:=\displaystyle{\sup_{0 \leq s \leq t}}  \lambda_{s}<\infty.
\end{split}	
\end{equation}
Furthermore,} we assume that $t\mapsto\sigma_{t}$ is continuous.
\end{assumption}
The following notation will be needed:
\begin{equation}\label{eq:process_average_notations}
	\overline{\sigma}^{2}_{t,h}:=\frac{1}{h}\int_{t}^{t+h}\sigma^{2}_{u}du, \quad 
	\overline{\gamma}_{t,h}:=\frac{1}{h}\int_{t}^{t+h}\gamma_{u}du,\quad
	\overline{\lambda}_{t,h}:=\frac{1}{h}\int_{t}^{t+h}\lambda_{u}du.
\end{equation}
Note that with these notations, our model assumptions imply that, for any $t,h\geq{}0$ and $k\in\mathbb{N}$,
$$
X^{c}_{t+h}-X^{c}_{t} 
=_{D} N \left( h\overline{\gamma}_{t,h}, h\overline{\sigma}^{2}_{t,h} \right),\qquad \Px \left( X_{t+h} - X_{t} \in dx | N_{t+h} - N_{t} = k \right) = \phi_{t,h}*f^{*k}(x)dx,
$$
where $\phi_{t,h}$ is the density of $X^{c}_{t+h}-X^{c}_{t}$, i.e.
$
\phi_{t,h}(x):= \frac{1}{\overline{\sigma}_{t,h}\sqrt{h}}\phi \left( \frac{x- h\overline{\gamma}_{t,h}}{\overline{\sigma}_{t,h}\sqrt{h}} \right)
$.
For these types of processes, the associated local characteristics are of the form $(\gamma,\sigma,\nu)$, where the density of the local L\'{e}vy measure is given by $\nu_{t}(x) = \lambda_{t} f(x)$. 

\begin{assumption}\label{assumption:min_f_zero}
{
The jump density $f$ {\Blue has} the form 
\begin{equation} \label{eq:jump_density_mixture_form}
f(x) = pf_{+}(x){\bf 1}_{[x \geq 0]} + qf_{-}(x){\bf 1}_{[x < 0]},
\end{equation}
where $p \in [0,1]$ and $q:=1-p$, and $f_{+}:[0,\infty)\to[0,\infty)$ and $f_{-}:(-\infty,0]\to[0,\infty)$ are \emph{bounded} functions such that 
$
\int_0^\infty f_+(x) dx
=
\int_{-\infty}^0 f_-(x) dx
=
1
$.
Furthermore, we assume that
\[
f_{\pm}(0)
=
\lim_{x\to{}0^{\pm}}f_{\pm}(x)
\in (0, \infty).
\]}
\end{assumption}
{The following notations will also be needed:}
\begin{equation}\label{eq:density_notation}
\mathcal{C}_{0}(f)
:= 
\lim_{\varepsilon \to 0^{+}} \frac{1}{2\varepsilon} \int_{-\varepsilon}^{\varepsilon} f(x) dx
=
p f_{+}(0) + q f_{-}(0),\quad
\mathcal{C}_d(f)
:=
\left| p f_{+}(0) - q f_{-}(0) \right|,
\quad
\mathcal{C}_m (f) 
:= 
\min \lbrace{f_{+}(0),f_{-}(0) \rbrace}.
\end{equation}
Note that $\mathcal{C}_{0}(f) = f(0)$ and $\mathcal{C}_d(f)=0$ if $f$ is continuous at the origin. 
For some results, we also need the following assumption:
\begin{assumption}\label{assumption:min_f_prime_zero}
$f_{+} \in {\rm C}^1([0, b))$, 
$f_{-} \in {\rm C}^1((a, 0])$,
for some {$a \in (-\infty, 0)$, $b \in (0, \infty)$}
and
$f_{\pm}'(0)
:= 
\lim_{x \to 0^{\pm}} f'_{\pm}(x)$
exists.
\end{assumption}

Throughout, we assume that we observe the process $X$ {\Blue at evenly spaced times,
\begin{equation}\label{DfnHn}
	t_{i}:=ih_{n},\quad i=0,\dots,n,
\end{equation}
where $h_{n}$ is the time span between observations and $T:=T_{n}:=nh_{n}$ is the time horizon. We will also} use {\Blue $\Delta^n_i X := X_{t_{i}} - X_{t_{i - 1}}$ to denote the increment of the underlying process over $[t_{i-1},t_{i})$}, and when no ambiguity can be brought, we will drop the superscript $n$. 
Finally, we introduce the {\Blue jump detection procedure we consider in this work. We first specify} a vector of {thresholds} $[\vec{B}]_{T}^{n} = (B^n_1, ..., B^n_n)$, where we often drop the superscript $n$ when no confusion can be generated. {\Blue Given} $[\vec{B}]_{T}^{n} $, we {\Blue would} conclude that a jump {\Blue had occurred} during $[t_{i - 1}, t_i)$ {\Blue whenever} $|\Delta_i X| > B_i$. {\Blue As a byproduct of this jump detection criterion, we can then devise the following natural} estimators {of $N_{T}$, $J_{T}$, and  the integrated variance $IV_T:=\int_0^T \sigma_s^2 ds$:}
\begin{equation}\label{eq:est_N_J_IV}
{\widehat{N}_T
= 
\sum_{i = 1}^n \textbf{1}_{\{|\Delta_i X| > B_i\}}, \quad
\widehat{J}_T
= 
\sum_{i = 1}^n (\Delta_i X) \textbf{1}_{\{|\Delta_i X| > B_i\}}, \quad
\widehat{IV}_T
= 
TRV(X)[\vec{B}]_{T}^{n}
=
\sum_{i = 1}^n (\Delta_i X)^2 \textbf{1}_{\{|\Delta_i X| \leq B_i\}}}.
\end{equation}
{These estimators were first studied in \cite{mancini2001disentangling}, \cite{mancini2004estimation}}.  The estimator $\widehat{IV}_T$ has extensively been studied in the literature and is commonly called the truncated or thresholded realized quadratic variation (TRV) of $X$. 

\subsection{Optimal Threshold and Its Approximation}\label{section:optimal_thresholding}

In this subsection, we {\Blue formulate} the problem of optimal threshold selection. {We adopt the approach in \cite{FLN2013}, which we now {briefly review} for completeness. We} seek to find a threshold $[\vec{B}]_{T} = (B_1, ..., B_n) \in \mathbb{R}_{+}^{n}$ to minimize the {loss function}:
\begin{equation}\label{eq:loss_function}	
L([\vec{B}]_{T})
:=
\Ex \left(
\sum_{i=1}^{n} \textbf{1}_{ \{ |\Delta_{i}X| > B_{i}, \Delta_{i}N= 0 \} } 
+ 
\sum_{i=1}^{n} \textbf{1}_{ \{ |\Delta_{i}X| \leq B_{i}, \Delta_{i}N\neq 0 \}}  \right).
\end{equation}
The above loss function {represents} the expected number of {``jump"} {mis-classifications} (i.e., subintervals erroneously classified as having {jumps when in fact they do not, or not having jumps when in fact they do)}.  The  previous formulation gives the same weight to both types of error, while a more general loss function is given by:
\begin{equation}\label{eq:loss_function_general}
L([\vec{B}]_{T}; w)
:=	
\Ex \left(
\sum_{i=1}^{n} \textbf{1}_{ \{ |\Delta_{i}X| > B_{i}, \Delta_{i}N= 0 \} } 
+ w
\sum_{i=1}^{n} \textbf{1}_{ \{ |\Delta_{i}X| \leq B_{i}, \Delta_{i}N\neq 0 \}}  \right).
\end{equation}
For our purpose, \eqref{eq:loss_function} is enough, but in {certain applications, \eqref{eq:loss_function_general} may} be useful. For instance, it is more likely that {market participants} become more conservative when they {erroneously} identify a price {change} as an unusual event, i.e., a jump. In this case, {one may prefer to} take $w <1$.

In both \eqref{eq:loss_function} and \eqref{eq:loss_function_general}, the loss function is additive. Therefore, we can optimize each $B_i$ separately. Indeed, we define the following loss function for given $t$ and $h$:
\begin{align}\label{eq:loss_single_term}
L_{t,h}(B; w) 
:= 
\Px \left( |X_{t+h}-X_{t}| > B , N_{t+h}-N_{t} = 0 \right) 
+ 
w \Px \left( |X_{t+h}-X_{t}| \leq B , N_{t+h}-N_{t} \neq 0 \right).
\end{align}
If we {were} able to devise a method to find $B^* = \mbox{argmin}_B L_{t,h}(B; w)$ for any $t$ and $h$, then, by setting $t = t_{i - 1}$ and $h = t_{i} - t_{i - 1}$, we {would be} able to {\Blue specify} the whole optimal $[\vec{B}]_T$.
Obviously, the first issue that we have to address is whether or not there is a global minimum point {\Blue $B^{*}$}.
As it turns out, the loss function \eqref{eq:loss_single_term} is {quasi-convex}\footnote{{A mapping $g: D \rightarrow \mathbb{R}$, for convex $D$,  is quasi-convex if for any $\lambda \in [0,1]$ and $x,y \in D$, $g(x \lambda + y (1-\lambda)) \leq \max \lbrace g(x), g(y) \rbrace$.}} in $B$, when $h$ is small enough. This property was established in  \cite{FLN2013} for a \emph{driftless} L\'evy processes (i.e., $\gamma\equiv{}0$ and $\sigma$ and $\lambda$ are constants). Nonzero drifts create some {nontrivial} subtleties that are resolved in the following {\Blue theorem, which was stated without proof in  \cite{FLN2013}.}
\begin{theorem}[\textbf{Uniform {Quasi-Convexity} of the Loss Functions}] \label{thm:uniform_convexity}
Assume that we have model \eqref{eq:jump_diffussion_model}, and
{Assumptions} {\ref{assumption:boundedness_of_intensity}-\ref{assumption:min_f_prime_zero}} are satisfied. {Then, for} any fixed $T > 0$, there exists $h_{0}:=h_{0}(T) > 0$, such that, for all $t \in [0,T]$, $h \in (0,h_{0}]$, and $w > 0$, the function $L_{t,h}(B;w)$ is \textbf{quasi-convex} in $B$, and possesses a unique global minimum point $B^{*}_{t,h}$. 
\end{theorem}
We proceed to give a fixed-point formulation of the optimal threshold $B^{*}_{t,h}$, which in turn {enables us to find a {second-order} asymptotic expansion for $B^{*}_{t,h}$ in a high-frequency asymptotic regime ($h\to{}0$).}  {This} characterization will {equip us with the} theoretical basis for developing feasible estimation algorithms later. {\Blue In what follows, we focus on the case of $w=1$ and for easiness of notation, we drop the variable $w$ in $L_{t,h}(B; w)$.}

\begin{theorem} [\textbf{Characterizations of the Optimal Threshold}] \label{thm:optimal_threshold_characterizations}
Assume that we have model \eqref{eq:jump_diffussion_model}, and
{Assumptions \ref{assumption:boundedness_of_intensity}-\ref{assumption:min_f_prime_zero}} are satisfied.
For each fixed $T>0$, there exists $h_{0}:=h_{0}(T)>0$ such that, for any $t\in [0,T]$ and $h\in (0,h_{0})$, the optimal threshold $B_{t,h}^{*}$, based on the increment $X_{t+h}-X_{t}$, is such that,
\begin{align}\nonumber
B_{t,h}^{*} &=  h\overline{\gamma}_{t,h} + \sqrt{2 h\overline{\sigma}^{2}_{t,h}} \left[\ln \left(1+\exp \left( \frac{-2B_{t,h}^{*}\overline{\gamma}_{t,h}}{\overline{\sigma}^{2}_{t,h}}\right)\right)\right.\\
&\qquad\qquad\qquad\quad\qquad \left.-\ln\left(\sqrt{2\pi h\overline{\sigma}^{2}_{t,h}}\sum_{k=1}^{\infty} \frac{\left( h\overline{\lambda}_{t,h}\right)^{k}}{k!} \left[\phi_{t,h}*f^{*k}(B_{t,h}^{*})+\phi_{t,h}*f^{*k}(-B_{t,h}^{*})\right] \right)\right]^{1/2}. \label{eq:fixed_point_equation}
\end{align}
Furthermore,  
as $h \to{} 0$, we have {the asymptotics:}
\begin{equation}  \label{eq:optimal_threshold_asym}
B_{t,h}^{*}
=
\sqrt{h}\, {\bar\sigma_{t,h}} 
\left[ 
3\log \left( 1/h \right)
-
2 \log \left( \sqrt{2\pi}\, \mathcal{C}_0(f) {\bar\sigma_{t,h} \bar\lambda_{t,h}} \right)
\right]^{1/2}
+ 
o(h^{\frac{1}{2} + \alpha} ),
\end{equation}
for \emph{any} $\alpha \in (0, 1/2)$. {\Blue If, furthermore, $\sigma_{t}^{2},\lambda_{t}\in C^{1}((0,T))$ and continuous on $[0,T]$, then} the asymptotics in \eqref{eq:optimal_threshold_asym} remains true if we replace $\bar\sigma_{t,h}$ and $\bar\lambda_{t,h}$ with $\sigma_t$ and $\lambda_t$, respectively.
\end{theorem}

\begin{remark}
	The last assertion of Theorem \ref{thm:optimal_threshold_characterizations} remains true if $t\to\sigma_{t}^{2}$ is H\"older continuous for any exponent $\chi\in (0,1/2)$.  In particular, this is the case for any volatility model driven by a Brownian motion (see \cite[Ch.V, Exercise 1.20]{RevuzYor}). Recall that if a function is H\"older continuous with exponent $\chi$, then it is H\"older continuous for any exponent $\chi'<1/2$. 
\end{remark}
{The previous result extends the first-order approximation $\sqrt{
3\sigma_{t,h}^2 h \log \left( 1/h \right)}$ of \cite{FLN2013}, whose  
remainder is just of order $O\left(
 h^{1/2} \log^{-1/2} \left( 1/h \right)
\right)$. 
However, with} the second order approximation, the remainder is $o(h^{1 - \epsilon})$ for any {$\epsilon \in (1/2, 1)$}. It is convenient to introduce the following notations for the first- and second-order optimal threshold approximations, respectively:
\begin{equation}\label{eq:B1_B2_for_t_h}
B_{t,h}^{*1}
=
{\bar{\sigma}_{t,h}}\left[ 
3 h \log \left( 1/h \right)
\right]^{1/2}, \quad
B_{t,h}^{*2}
=
\sqrt{h}\, \bar{\sigma}_{t,h}
\left[ 
3\log \left( 1/h \right)
-
2 \log \left( \sqrt{2\pi} \mathcal{C}_0(f) \bar{\sigma}_{t,h} \bar{\lambda}_{t,h} \right)
\right]^{1/2}.
\end{equation}
{\Blue These tell us that, in a high-frequency sampling setting, the single most important parameter to determine a suitable threshold level $B$ is the spot volatility $\sigma_{t}$,} followed by the parameter  
$\nu_{t}(0):=\lambda_{t}\mathcal{C}_{0}(f)$, which broadly 
 determines the likelihood of a small jump occurrence around time $t$.
It is interesting to note that the optimal threshold ${B}^{*2}_{t,h}$ can differ substantially from $B^{*1}_{t, h}$ when $\sigma_{t} \lambda_{t}{\mathcal{C}_{0}(f)}$ is large.  
This is intuitive since, for instance, if $\sigma_{t}$ and ${C_{0}(f)}$ are fixed, as the jump rate $\lambda_{t}$ increases, the optimal threshold {\Blue should decrease} in order to account for an {\Blue increase in the appearance of ``small" jumps. If the threshold is not adjusted, there would be  more ``false-negatives", i.e., missed jumps.  On the other hand, as $\lambda_{t}$ decreases, the optimal threshold should be larger} in order to offset an increment in the likelihood of {false-positives} (namely, wrongly concluding the occurrence of a jump during the small interva $[t,t+h]$). {\Red Similarly, for fixed $\sigma_{t}$ and  $\lambda_{t}$ as the likelihood for small jumps, approximately parameterized by ${C_{0}(f)}$, increases (decreases) the optimal threshold decreases (increases) accordingly.}

Although we have proved the asymptotic properties of \eqref{eq:B1_B2_for_t_h}, these optimal thresholds are not yet {feasible,} since we still need to estimate the spot volatility $\sigma_t^2$, {\Red jump} intensity $\lambda_t$, and the {mass concentration of the jump density at the origin, $\mathcal{C}_0(f)$}. We will introduce estimators to these quantities in {Subsection \ref{Sect:TKJD} and Section \ref{sec:threshold_kernel_spot_vol}}, {respectively.}

\begin{remark}
Although {the criterion (\ref{eq:loss_single_term})} provides a reasonable approach {for threshold selection,} there is no guarantee that {the resulting} optimal threshold is the one that minimizes the {mean-square error of the truncated realized quadratic variation $\widehat{IV}_{T}$ introduced in (\ref{eq:est_N_J_IV}).}  We refer to \cite{FLM2017} for some results regarding {the latter problem.}
\end{remark}

\subsection{{Bias and Variance}}\label{BiasVarianceSect}
We conclude with the following asymptotic result of the estimation error of the TRV, which generalizes {a result of} \cite{FLN2016} to non-homogeneous drift, volatility, and jump intensities. {\Blue As usual, the notation $a_{h}\sim b_{h}$, as $h\to{}0$, means that $\lim_{h\to{}0}a_{h}/b_{h}=1$.}
\begin{proposition} \label{prop:iv_bias_variance} 
Suppose that the assumptions of Theorem \ref{thm:optimal_threshold_characterizations} are enforced {and that  $\vec{B}=(B_{n})_{n\geq{}1}$ is set to be $B_{n}^{*1} = \sqrt{3 {{\sigma}_{t_i}^2} h_n \log(1/h_n)}$.
Then, as $n\to\infty$,}
\begin{align*}
\Ex \left[TRV(X)[\vec{B}]_{T}^{n} \right] - \int_0^T \sigma_s^2 ds  
\sim 
h_{n} \int_0^T (\gamma_s^2 - \lambda_s \sigma_s^{2}) ds, 
\quad 
{\rm Var}\left( TRV(X)[\vec{B}]_{T}^{n} \right) 
\sim 
2 h_{n} \int_0^T\sigma_s^{4} ds. 
\end{align*}
Furthermore, the asymptotic {\Blue behavior} above {also} holds with any {threshold sequence} of the form 
\[
\check{B}_{n,i} = \sqrt{{c_{n,i} {\sigma}_{t_{i}}^2} {h_{n}} \log(1/{h_{n}})} + o(\sqrt{{h_{n}} \log(1/{h_{n}})}),
\] provided that {$c:=\liminf_{n\to\infty}\inf_{i}c_{n,i} \in (2, \infty)$}.
\end{proposition}
\begin{proof}
{Let us write $B_{n,i}$ of the form  $\sqrt{3 {{\sigma}^{2}_{{t_{i}}}} h_n \log(1/h_n)}$}. {The bias of the} TRV estimator can be decomposed as the following:
\begin{align}\nonumber
&\quad
TRV(X)[\vec{B}]_{T_{n}}^{n} - \int_0^T \sigma_s^2 ds \\
& =
\sum_{i = 1}^n \left(
 | \Delta^n_i X |^2 \textbf{1}_{[\Delta^n_i N = 0]} - h_n \overline{\sigma}^{2}_{t_{i - 1},h_n}
\right)
+
\sum_{i = 1}^n | \Delta^n_i X |^2 \textbf{1}_{[ | \Delta^n_i X | \leq B_n, \Delta^n_i N \neq 0]}
-
\sum_{i = 1}^n | \Delta^n_i X |^2 \textbf{1}_{[ | \Delta^n_i X | > B_n, \Delta^n_i N = 0]}.\label{ScnRel}
\end{align}
{\Blue Using Lemmas C.1 and C.2 in \cite{FLN2016}} as well as Assumption \ref{assumption:boundedness_of_intensity}, for any $0<\epsilon<1/2$, we have:
\begin{align}\label{T1In}
\mathbb{E}  \left[ 
| \Delta^n_i X |^2 \mathbf{1}_{[ | \Delta^n_i X | \leq B_{n,i}, \Delta^n_i N \neq 0]}
\right] 
& = O(B_{n,i}^{3} h_{n})=
O(h_{n}^{\frac{5}{2}} 
\left[ 
\log \left( 1/h_{n} \right)
\right]^{\frac{3}{2}}), \\
\label{T2In}
\mathbb{E} \left[ 
|\Delta^n_i X |^2 \mathbf{1}_{[ | \Delta^n_i X | > B_{n,i}, \Delta^n_i N = 0]}
\right]
& = O(\sqrt{h_{n}}B_{n,i}\phi(B_{n,i}/\bar{\sigma}_{t_{i},h_{n}}\sqrt{h_{n}}))=
O \left( h_{n}^{\frac{5}{2}-{\epsilon}} 
\left[ 
\log \left( 1/h_{n} \right)
\right]^{\frac{1}{2}} \right),
\end{align}
where the $O(\cdot)$ terms are uniform in $i$. These would imply {that the second and third terms of (\ref{ScnRel}) are of orders $O_P(h^{3/2}
\left[ 
\log \left( 1/h \right)
\right]^{3/2})$ and $O_P( h^{3/2-{\epsilon}} 
\left[ 
\log \left( 1/h \right)
\right]^{1/2} )$, respectively. Both of these terms are then $o(h_{n})$.}
For the first term therein, note that
\begin{align*}
\Ex \left[
 | \Delta^n_i X |^2 \mathbf{1}_{[\Delta^n_i N = 0]} - h_n \overline{\sigma}^{2}_{t_{i - 1},h_n}
\right]
& =
\Px (\Delta^n_i N \neq 0) h_n \overline{\sigma}^{2}_{t_{i - 1},h_n} 
+
\Px (\Delta^n_i N = 0) h_n^2 \overline{\gamma}^{2}_{t_{i - 1},h_n} \\
& = 
h_n^2 (\overline{\gamma}^{2}_{t_{i - 1},h_n} - \overline{\lambda}_{t_{i - 1},h_n} \overline{\sigma}^{2}_{t_{i - 1},h_n} ) + O(h_n^3) , \\
\Ex \left[ \left(
 | \Delta^n_i X |^2 \mathbf{1}_{[\Delta^n_i N = 0]} - h_n \overline{\sigma}^{2}_{t_{i - 1},h_n}
\right)^2 \right]
& =
\Px (\Delta^n_i N \neq 0) h_n^2 \overline{\sigma}^{4}_{t_{i - 1},h_n} 
+
\Px (\Delta^n_i N = 0) 2 h_n^2 \overline{\sigma}^{4}_{t_{i - 1},h_n} 
\\&=
2 h_n^2 \overline{\sigma}^{4}_{t_{i - 1},h_n} + O(h_n^3).
\end{align*}
Calculating the summation of the above and noticing the independence of different terms, we conclude {the} first part of the desired result.
{For {$\check{B}_{n,i}$, the term (\ref{T2In}) will instead be of order 
$O_P( h^{1+c/2-{\epsilon}} 
\left[ 
\log \left( 1/h \right)
\right]^{1/2})$.} 
Therefore,} as long as {$c > 2$}, the asymptotic behavior does not change. This proves the second part of the desired result.
\end{proof}

\begin{remark}
The motivation for considering the threshold {$\breve{B}_{n,i}$} in {Proposition \ref{prop:iv_bias_variance}} {comes from} the fact that the true value {of $\sigma^{2}$} is not available and, in practice, we have to use an estimate $\hat{\sigma}^{2}$ of it. Suppose we have an estimator of $\sigma_{t_i}^2$ denoted by $\hat{\sigma}_{t_i}^2$, and we use the corresponding estimated threshold $\hat{B}_{n}^{*1} = \sqrt{3 \hat{\sigma}_{t_i}^2 h_n \log(1/h_n)}$. 
The second part of Proposition \ref{prop:iv_bias_variance} tells us that if, {for instance, the estimator is such that} $\liminf_{n \to \infty} \hat{\sigma}_{t_{i}}^2 / \sigma_{t_{i}}^2 = c > 2/3$, we would have 
$
\hat{B}_{n}^{*1} = \sqrt{3 \hat{\sigma}_{t_i}^2 h_n \log(1/h_n)}
\geq
\sqrt{3 (c - \epsilon) \sigma_{t_i}^2 h_n \log(1/h_n)},
$
for $n$ large enough and $\epsilon \in (0, c - 2/3)$.
This {will result} in an estimator such that the asymptotics of {the} expectation and variance of Proposition \ref{prop:iv_bias_variance} hold.
\end{remark}

\subsection{A Threshold-Kernel Estimation of the Jump Density at $0$}\label{sec:jump_density_estimation}\label{Sect:TKJD}

In this section, we investigate the estimation of the jump density at the origin, which is needed {in order to implement the second order optimal threshold $B^{*2}_{t, h}$ given by \eqref{eq:B1_B2_for_t_h}}. We propose a method based on kernel estimators. {\Blue For a related method, but for a more general class of It\^o semimartingales, see \cite{Ueltzhofer}. The main difference between the method proposed below and the one proposed in that paper is the thresholding technique.}

{\Blue We impose} the following regularity conditions, which in particular imply that $\mathcal{C}_{0}(f)=f(0)$. 
\begin{assumption}\label{assumption:smooth_density}
$f \in C^2 \left( [a, b] \right)$ for {some} $a < 0 < b$. Also, $f(0) \neq 0$ and $f^{\prime\prime} (0) \neq 0$.
\end{assumption}

\begin{remark}
{It is possible to relax the previous assumption. {For instance, if the density  $f$ merely satisfies} Assumption \ref{assumption:min_f_zero}, the} estimation of $f(0^+)$ and $f(0^-)$ would have to be done separately using one-sided kernel estimators. The basic idea is the same {as} what we present below, but the convergence rate and the choice of bandwidth will be different.
\end{remark}
{As mentioned above, we wish to construct a consistent estimator for $\mathcal{C}_0(f)=f(0)$, which is not feasible during a 
fixed time interval $[0, T]$. 
Hence, in this part, 
we consider a high-frequency/long-run sampling setting, where {simultaneously} 
\[
	{h}_{n}=t_{i}-t_{i-1} \to 0,\qquad T_{n}=t_{n} \to \infty,
\]
as $n\to\infty$. Throughout, we also assume that $\gamma$, $\sigma$, and $\lambda$ are constant so that the distribution of $\Delta_i X$ does not depend on $i$.

In the spirit of threshold {estimation, the basic idea is to treat the {``large"} increments $\Delta_i X$, whose absolute values exceed an appropriate threshold, as {proxies of} the process' jumps. These large increments can then be plugged into a {standard} kernel estimator of $ f(0)$.  {Concretely, we consider the
estimator:}
\begin{equation}\label{eq:threshold_kernel_est_jump_density}
{2\hat{f}(0)}:=\frac{1}{|\{{i}: |\Delta_i X | > B\}|} 
\sum_{\{{i:} | \Delta_i X | > B \} } K_{\delta}( |\Delta_i X| - B) ,
\end{equation}
under the convention that $0/0=0$ in the case that $\{{i:} | \Delta_i X | > B \} =\emptyset$. As usual, {$K_{\delta}(x):= K(x/{\delta})/{\delta}$, where {$K:[0,\infty)\to [0,\infty)$} is a right-sided kernel function such that $\int_0^\infty K(x) dx = 1$ and ${\delta}$} is the bandwidth parameter. We also use $|A|$ to denote the number of elements in {a} set $A$. 
We expect that the estimator \eqref{eq:threshold_kernel_est_jump_density} will have poor performance if $|\{i: |\Delta X_{i}| > B\}|$ 
is small, but, since we assume that $T \to \infty$ and $f(x)\neq{}0$ in a neighborhood of $\{0\}$, {for large-enough $n$, we have $P(\{|\Delta_{i} X| > B\} = \emptyset) \approx e^{-\lambda T} \to 0$. For} our implementation of (\ref{eq:threshold_kernel_est_jump_density}) {in the Monte Carlo studies of Section \ref{sec:monte_carlo_study}}, 
we will set $\hat{f}(0)=0$ if $|\{{i}: |\Delta_i X | > B\}|\leq{}5${, which simply makes the second order threshold to be the first order threshold}.

In what follows, $f^*$ {stands} for the density of {$|\Delta_{i} X|$}, {which depends on $n$}, while $f^*_{|\Delta X| | |\Delta X| > B}$ {stands} for the density of {$|\Delta_{i} X|$} conditioning on $|{\Delta_{i} X}| > B$.
{To analyze the performance of the estimator (\ref{eq:threshold_kernel_est_jump_density}) and choose a suitable thresholding level $B$ {and bandwidth ${\delta}$}, we decompose the estimation error into the following two terms:}
\begin{enumerate}[(i)]
\item
$E_1 = \frac{1}{|\{ |\Delta_{i} X| > B\}|} \sum_{ | \Delta_{i} X | > B} K_{{\delta}}( |\Delta_{i} X| - B) - f^*_{ |\Delta X| | |\Delta X| > B }(B)$,
\item
$E_2 = f^*_{ |\Delta X| | |\Delta X| > B }(B) - 2f(0)$.
\end{enumerate}

{Next, we {follow} a ``greedy" strategy to {determine suitable values for the} threshold $B$ and bandwidth ${\delta}$. Specifically, we minimize $E_2$ to obtain an ``optimal" threshold $B$, and with that given, we minimize $E_1$ to obtain {an} ``optimal" bandwidth ${\delta}$. Minimizing $E_1 + E_2$ directly will be a much more involved problem, and requires more assumptions. However, we believe solving such a problem does not significantly improve the performance of the proposed estimator. Therefore, we leave it as an open problem.

{Minimizing $E_1$ over ${\delta}$} \emph{given $B$} is {closely related to the standard} theory of kernel density estimation, so we can directly {apply the general theory for such a problem}. We only need to ensure that $|\{ |\Delta_{i} X| > B\}| \to \infty$, which follows from Proposition \ref{prop:conditions_for_convergence_E2} {\Blue below} with {the} additional assumption that $T\to\infty$.
{Two} widely used methods are plug-in method and cross-validation, which both have pros and cons. These methods are beyond the scope of this paper {and, for simplicity, we instead use the well-known Silverman's (1986) rule of thumb for bandwidth selection}:
\begin{equation}\label{eq:jump_density_bandwidth_rule_of_thumb}
{\delta}
=
1.06 L^{-1/5} \mbox{sd},
\end{equation}
where ``$ \mbox{sd}$" is the standard deviation of $\{\Delta_i X : |\Delta_i X| > B \}$ 
and $L$ is the number of observations, i.e. $| \{\Delta_i X : |\Delta_i X| > B \} |$. Such a rule of thumb works the best with Gaussian kernel function and Gaussian density function. However, the method {is known to be} robust for other kernel and density functions.

We now proceed to show that $B^{*} = \sqrt{4 {h} \sigma^2 \log(1/{h})}$ minimizes the leading order terms of the second error $E_2$.  The proof of the following two results are given in Appendix \ref{sec:appendix_optimal_thresholding}. 
 \begin{proposition}\label{prop:conditions_for_convergence_E2}
Suppose that Assumption \ref{assumption:smooth_density} is satisfied and $\gamma$, $\sigma$, and $\lambda$ are constant. {Further assume that $B \to 0$ and $B/\sqrt{{h}} \to \infty$.} Then, $E_2$ converges to $0$ as ${h}\to{}0$ if and only if ${h}^{-3/2} \exp\left(  -\frac{B^2}{2 {h} \sigma^2 } \right) \to 0$. {Under this condition, we have}
\begin{equation}\label{eq:kernel_jump_density_E2}
E_2
=
\frac{ 2 }{\lambda \sqrt{2 \pi {h}^3 \sigma^2}} \exp\left(  -\frac{B^2}{2 {h} \sigma^2 }\right)
+
2f(0) B + o(B)+o(h^{-3/2}e^{-\frac{B^{2}}{2h\sigma^{2}}}).
\end{equation}
Furthermore, if $E_2$ converges to $0$, then $\Px ( |{\Delta_{i} X}| > B ) = \lambda {h} +  o({h})$, {as ${h}\to{}0$.}
\end{proposition}

In addition to {providing us conditions} for the error $E_{2}$ to vanish,
Proposition \ref{prop:conditions_for_convergence_E2} implies that, {in that case,} 
$\mathbb{E} \left[ |\{ i : |\Delta_i X| > B\}| \right] = \lambda T + o(T)$, {as $T\to\infty$ and ${h}\to{}0$}. Therefore, the average sample size that can be used for the estimation of $f(0)$ is approximately constant with respect to $B$. {Heuristically, this suggests that the selection of $B$ will not affect significantly the selection of {\Blue $\delta$} that minimizes $E_{1}$.} We {are now ready to} obtain an approximate optimal threshold $B$, which minimizes the leading order terms of $E_2$. 

\begin{corollary}\label{lemma:approximate_optimal_threshold_E2}
The approximate optimal threshold {\Blue $\widetilde{B}^{*}$} that minimizes the leading order term of $E_2$ given by \eqref{eq:kernel_jump_density_E2} is {such that}
\begin{equation}\label{eq:approximate_optimal_threshold_E2}
{\Blue \widetilde{B}^*} 
= 
\sqrt{4 {h} \sigma^2 \log(1/{h})} + O(\sqrt{{h} \log\log(1/{h})}),
\end{equation}
\end{corollary}

It is interesting to notice that the ``optimal" threshold here is not the same as the one identified in the previous section. Indeed, if we do use the optimal threshold $B^{*1}$ or $B^{*2}$ in \eqref{eq:B1_B2_for_t_h}, $E_2$ would {diverge}. 
It is interesting and important to get some sense why the optimal thresholds differ from each other. Indeed, in the previous section, we optimize the expected number of jump misclassification. In that case, we are minimizing the sum of unconditional false positive (mistakenly claim a jump) and unconditional false negative (miss a jump). However, since the probability that a jump occurs is so small, proportional to the length of the time increments, the probability of having a false negative, {by nature}, cannot be too large. Therefore, by having the expected number of misclassification as the objective function, we would choose a threshold in favour of having a much smaller unconditional false positive rate. 
{As it turns out, 
if we choose $B^{*1}$ or $B^{*2}$, conditioning on $|\Delta X| > B$, the probability} that no jump occur is comparable to the probability that a jump occurs, both $O({h})$. That is, the conditional false negative rate does not vanish. Such a situation would minimize the expected number of {\Red misclassifications}, but would not enable us to distinguish the distribution of the jump from the noise. Using $\sqrt{4 {h} \sigma^2 \log(1 / {h})}$, on the other hand, makes the conditional false negative {vanishing and, thus, enables us to get consistent} estimation of jump density.}

\section{Threshold-Kernel Estimation of Spot Volatility}\label{sec:threshold_kernel_spot_vol}

{In this section, we consider the estimation of the spot volatility of a jump-diffusion process, which is needed to implement the approximate optimal threshold formula (\ref{eq:optimal_threshold_asym}), but is also an important problem on its own. Unlike Section \ref{sec:optimal_threshold}, here we {\Blue also work with certain} stochastic volatility models. The precise conditions are given below.}  
 
The {idea} of kernel estimation of spot volatility is to take a weighted average of the squared increments (see, e.g., \cite{foster1994continuous} and \cite{fan2008spot}):
\begin{equation}\label{eq:original_kernel_spot_vol}
{\hat{\sigma}_{\tau}^{2}}:=KW(\tau, n, {\delta})
:=
\sum _{i=1}^n K_{{\delta}}(t_{i-1} - \tau) (\Delta_i X)^2.
\end{equation}
{Here,} $K(\cdot)$ is {a} kernel function with $\int K(x) dx = 1$, $K_{\delta}(x)=K(x/\delta)/\delta$, and ${\delta} > 0$ is the bandwidth.
However, when jumps do occur, the estimator above becomes inaccurate. A natural idea is to combine \eqref{eq:original_kernel_spot_vol} with the threshold method.
Concretely, {given a} threshold vector $[\vec{B}]^{n}_{T}=(B_{1}^{n},\dots,B_{n}^{n})$, we  consider the local threshold-kernel estimator:
\begin{equation}\label{eq:threshold_kernel_spot_vol}
{\hat{\sigma}_{\tau}^{2}}:=TKW(\tau, n, {\delta})
:=
\sum _{i=1}^n K_{{\delta}}(t_{i-1} - \tau) (\Delta_i X)^2 \textbf{1}_{\{|\Delta_i X| \leq B_i \}}.
\end{equation}
{In what follows} {we will} investigate the properties of \eqref{eq:threshold_kernel_spot_vol}. In order to do this, we will have to deal with the randomness of the volatility, {for which we extend some of the results} in  \cite{FLL2016}. 
We will mention the assumptions on $\{\sigma_t\}_{t \geq 0}$ and $K$ in Subsection \ref{sec:assumption_stochastic_vol}, and then discuss the asymptotic properties of \eqref{eq:threshold_kernel_spot_vol} in subsequent subsections.

\subsection{Assumptions on the Volatility Process}\label{sec:assumption_stochastic_vol}

The first {\Blue assumptions are some non-leverage and boundedness conditions, which enable us to condition on the whole path of the volatility and drift and use estimates from \cite{FLN2016}}:
\begin{assumption}\label{IndependentBoundednessCondition}
In \eqref{eq:jump_diffussion_model}, {$(\gamma, \sigma)$ are locally bounded c\'adl\'ag} independent of the Brownian motion $W$ and {\Blue the jump component $J$. Furthermore, there exists a deterministic $M_T <\infty$ for which $\bar{\gamma}_{T}$ and $\bar{\sigma}_{T}^{2}$ defined in (\ref{BdnessOnsigmaetc}) satisfy $\bar{\gamma}_{T}<M_{T}$ and $\bar{\sigma}_{T}^{2}<M_{T}$. 
The intensity $\lambda$ is still assumed to be deterministic such that $\underline{\lambda}_{t}:=\inf_{0 \leq s \leq t}  \lambda_{s}> 0$ and $\bar{\lambda}_{t}:=\sup_{0 \leq s \leq t}  \lambda_{s}<\infty$.}
\end{assumption}
{We now introduce the} key assumption on the volatility process.
\begin{assumption}\label{VolatilityConditionGeneral}
Suppose that for $\varpi > 0$ and certain functions $L: \mathbb{R}_+ \rightarrow \mathbb{R}_+$, $C_{\varpi}: \mathbb{R} \times \mathbb{R} \rightarrow \mathbb{R}$, such that $C_{\varpi}$ is not identically zero and
\begin{equation}\label{HomoFunction}
\begin{split}
C_{\varpi}(hr,hs) & = h^{\varpi}C_{\varpi}(r,s), \quad
\mbox{ for } r, s \in \mathbb{R}, h \in \mathbb{R}_+,
\end{split}
\end{equation}
the variance process $V := \{V_t = \sigma_t^2 : t\geq 0\}$ satisfies
\begin{equation}\label{SCDefinitionGeneral}
\mathbb{E} [ (V_{t+r} - V_{t})(V_{t+s} - V_{t}) ] 
= L(t)C_{\varpi}(r,s) + o ((r^2 + s^2)^{\varpi/2}) , \quad r,s \rightarrow 0.
\end{equation}
\end{assumption}

An additional assumption on the kernel function $K$ is the following:
\begin{assumption}\label{AdmissibleKernel}
Given $\varpi > 0$ and $C_{\varpi}$ as defined in Assumption \ref{VolatilityConditionGeneral}, {the kernel} function $K:\mathbb{R} \rightarrow \mathbb{R}$ satisfies the following conditions:
\begin{enumerate}[(1)]
\item
$\int K(x) dx = 1${;}
\item
$K$ is Lipschitz and piecewise $C^1$ on its support $(A,B)$, where $-\infty \leq A < 0 < B \leq \infty${;}
\item
{(i)} $\int |K(x)||x|^{\varpi} dx < \infty$; {(ii) $K(x)x^{\varpi + 1} \rightarrow 0$, as} $|x| \rightarrow \infty$; {(iii)} $\int |K^{\prime}(x)| dx < \infty$, {(iv)} $V_{-\infty}^{\infty} (|K^{\prime}|) < \infty$, where $V_{-\infty}^{\infty}(\cdot)$ is the total variation{;}
\item
$\iint K(x)K(y)C_{\varpi}(x,y)dxdy > 0$.
\end{enumerate}
\end{assumption}

We refer to \cite{FLL2016} for more details on {\Blue Assumptions} \ref{IndependentBoundednessCondition}, \ref{VolatilityConditionGeneral} and \ref{AdmissibleKernel}. {We just mention here that  Assumption \ref{VolatilityConditionGeneral} covers a wide range of} frameworks such as deterministic and smooth volatility, Brownian motion and fractional Brownian motion driven volatility, etc. In the following subsection, we will establish asymptotic properties of \eqref{eq:threshold_kernel_spot_vol} based on Assumption \ref{IndependentBoundednessCondition}, \ref{VolatilityConditionGeneral} and \ref{AdmissibleKernel}.

\subsection{Asymptotic Properties of Threshold-Kernel Estimator}\label{sec:asymptotic_threshold_kernel_est}

\cite{FLL2016} proves the following result under Assumption \ref{IndependentBoundednessCondition}, \ref{VolatilityConditionGeneral}, and \ref{AdmissibleKernel} (c.f. Section 3 therein):
\begin{equation}\label{eq:original_kernel_error}
\begin{split}
&\mathbb{E} \left[ \left(
\sum _{i=1}^n K_{{\delta}}(t_{i-1} - \tau) (\Delta_i X^c)^2
-
\sigma_\tau^2
\right)^2 \right] \\
&\quad
=
2\frac{{h}}{{\delta}} \mathbb{E}[\sigma_{\tau}^4] \int K^2(x)dx
 +
{\delta}^{\varpi} L(\tau) \iint K(x)K(y) C_{\varpi} (x,y) dxdy 
+
{o\left(\frac{{h}}{{\delta}}\right)} + {o \left({\delta}^{\varpi}\right)},
\end{split}
\end{equation}
where $X^c$ is the continuous part of $X$ {defined in (\ref{eq:jump_diffussion_model})}. The key result to extend the {theory of kernel estimators, as developed in \cite{FLL2016},} to the threshold-kernel estimators (\ref{eq:threshold_kernel_spot_vol}) is the following.
\begin{proposition}\label{thm:diff_of_kw_and_tkw}
Suppose that {\Blue Assumptions 
\ref{assumption:min_f_zero},} \ref{IndependentBoundednessCondition}, \ref{VolatilityConditionGeneral}, and \ref{AdmissibleKernel} are satisfied, {and take a bandwidth sequence ${\delta}_{n}$ such that ${h}_{n}/{\delta}_{n}\to{}0$. Let 
${B_{i}}:=B_{n,i}(c)
 := \sqrt{c {\bar{\sigma}_{t_{i},h}^2} {{h}} \log(1/{{h}})} + o({\sqrt{{h} \log(1/{h})}})
$, with $c > 0$.} Then, we have:
\begin{equation}\label{eq:diff_of_kw_and_tkw_c}
\bar{\mathcal{E}}_{n}:=\sum _{i=1}^n K_{{\delta}}(t_{i-1} - \tau) \left[(\Delta_i X^c)^2- (\Delta_i X)^2 \textbf{1}_{\{|\Delta_i X| \leq B_i \}}\right]
=
O_P\left( \max\{ {h}, {h}^{c/2} \log^{{1/2}}(1 / {h}) \} \right).
\end{equation}
{Furthermore,} 
\begin{align}\label{eq:diff_of_kw_and_tkw_cForL2}
{\mathbb{E}\left(\bar{\mathcal{E}}_{n}^{2}\right)=
	{O}\left(\frac{{h}^2}{{\delta}}\right)
+ {O}\left( \frac{{h}^{1 + \frac{c}{2}}}{{\delta}} [\log(1/{{h}})]^{\frac{3}{2}} \right)
+ {O}\left( {h}^{c} \log(1/{{h}}) \right).}
\end{align}
\end{proposition}

\begin{proof}
{Let $\mathcal{E}_{i}:=(\Delta_i X)^2 \textbf{1}_{\{|\Delta_i X| \leq B_i \}} - (\Delta_i X^c)^2$ and} observe that
\begin{equation}
\begin{split}
{\mathcal{E}_{i}}
& =
- ( \Delta_i X^c )^2 \mathbf{1}_{[ \Delta_i N \neq 0]} 
+
( \Delta_i X )^2 \mathbf{1}_{[ | \Delta_i X | \leq {B_{i}}, \Delta_i N \neq 0]}
-
( \Delta_i X )^2 \mathbf{1}_{[ | \Delta_i X | > {B_{i}}, \Delta_i N = 0]}
=:\mathcal{E}_{i,1}+\mathcal{E}_{i,2}+\mathcal{E}_{i,3}.
\end{split}
\end{equation}
Now, {\Blue conditioning on the paths of $\sigma$ and $\gamma$ and applying Lemmas C.1-C.2 in  \cite{FLN2016}}, the following holds:
\begin{equation}
\begin{split}
& {( \Delta_i X )^2 \mathbf{1}_{[ | \Delta_i X | \leq B_{i}, \Delta_i N \neq 0]} =  {O_P\left( B_{i}^{3}{h}\right)}=O_P\left( {h}^{5/2} [\log(1/{{h}})]^{3/2} \right)}, \\
& {( \Delta_i X )^2 \mathbf{1}_{[ | \Delta_i X | > B_{i}, \Delta_i N = 0]} = 
{\Blue O_{P}(\sqrt{h_{n}}B_{n,i}\phi(B_{n,i}/\bar{\sigma}_{t_{i},h_{n}}\sqrt{h_{n}}))}=
O_P\left( {h}^{1 + c/2} [\log(1/{{h}})]^{1/2} \right)},\\
&{ ( \Delta_i X^c )^2 \mathbf{1}_{[ \Delta_i N \neq 0]}  = O_{P}({h}^2)}.
\end{split}
\end{equation}
From Assumption \ref{IndependentBoundednessCondition}, the above holds uniformly over $1 \leq i \leq n$. Therefore, by Assumption \ref{AdmissibleKernel}, we have:
$$
\sum _{i=1}^n K_{{\delta}}(t_{i-1} - \tau) \left[ 
(\Delta_i X)^2 \textbf{1}_{\{|\Delta_i X| \leq B_i \}} 
-
(\Delta_i X^c)^2
\right]
=
O_P\left( \max\{ {h}, {h}^{c/2} \log^{1/2}(1 / {h}) \} \right).
$$
{For the second assertion of the theorem, first note that 
\begin{align*}
\mathbb{E}\left(\bar{\mathcal{E}}_{n}\right)&=\sum _{i=1}^n K_{{\delta}}(t_{i-1} - \tau) \mathbb{E}\left[\mathcal{E}_{i,1}+\mathcal{E}_{i,2}+\mathcal{E}_{i,3}\right]\\
&=\sum _{i=1}^n K_{{\delta}}(t_{i-1} - \tau) \left[{O}({h}^2)+
{O}\left( {h}^{\frac{5}{2}} [\log(1/{{h}})]^{\frac{3}{2}} \right)
+ {O}\left( {h}^{1 + \frac{c}{2}} [\log(1/{{h}})]^{\frac{1}{2}} \right)\right]\\
&= {O}\left({h}\right) + 
{O} \left( {h}^{\frac{c}{2}} [\log(1/{{h}})]^{\frac{1}{2}} \right).
\end{align*}
Similarly,
\begin{align*}
{\rm Var}\left(\bar{\mathcal{E}}_{n}\right)&=\sum _{i=1}^n K^{2}_{{\delta}}(t_{i-1} - \tau) {\rm Var}\left((\Delta_i X^c)^2
- (\Delta_i X)^2 \textbf{1}_{\{|\Delta_i X| \leq B_i \}}\right)\\
&\leq{}4\sum _{i=1}^n K^{2}_{{\delta}}(t_{i-1} - \tau) \left[\mathbb{E}(\mathcal{E}^{2}_{i,1})+\mathbb{E}(\mathcal{E}^{2}_{i,2})+\mathbb{E}(\mathcal{E}^{2}_{i,3})\right]\\
&=\sum _{i=1}^n K^{2}_{{\delta}}(t_{i-1} - \tau) \left[{O}({h}^3)
+{O}\left( {h}^{2 + \frac{c}{2}} [\log(1/{{h}})]^{\frac{3}{2}} \right)\right]\\
&= {O}\left(\frac{{h}^2}{{\delta}}\right)
+{O}\left( \frac{{h}^{1 + \frac{c}{2}}}{{\delta}} [\log(1/{{h}})]^{\frac{3}{2}} \right).
\end{align*}
We then conclude 
the result.}
\end{proof}

With {Proposition \ref{thm:diff_of_kw_and_tkw}}, {\Blue we get}  the following proposition, which characterizes the leading order terms of the MSE of the threshold-kernel estimator (\ref{eq:threshold_kernel_spot_vol}). This allows us to {perform} bandwidth and kernel function selection.
\begin{proposition}\label{prop:MSE_tkw}
Assume that {Assumptions} {\Blue 
\ref{assumption:min_f_zero},} \ref{IndependentBoundednessCondition}, \ref{VolatilityConditionGeneral}, and \ref{AdmissibleKernel} are satisfied, and take the threshold vector to be $B_{{n,i}}(c) = \sqrt{c \bar{\sigma}_{t_i,h}^2 {h} \log(1/{h})} + o(\sqrt{{h} \log(1 / {h})})$ for any $c\in (\frac{\varpi}{\varpi + 1}, \infty)$. 
Then, we have that, for each $\tau \in (0,T)$,
\begin{equation}\label{eq:MSE_tkw}
\mathbb{E} \left[ \left(
TKW(\tau, n, {\delta})
-
\sigma_\tau^2
\right)^2 \right]
=
2\frac{{h}}{{\delta}} \mathbb{E}[\sigma_{\tau}^4] \int K^2(x)dx
 +
{\delta}^{\varpi} L(\tau) \iint K(x)K(y) C_{\varpi} (x,y) dxdy 
+
o\left(\frac{{h}}{{\delta}}\right) + o \left({\delta}^{\varpi}\right).
\end{equation}
\end{proposition}

\begin{proof}
We consider the following decomposition:
\begin{equation}
\begin{split}
TKW(\tau, n, {\delta}) - \sigma_\tau^2
& =
\sum _{i=1}^n K_{{\delta}}(t_{i-1} - \tau) \left[ 
(\Delta_i X)^2 \textbf{1}_{\{|\Delta_i X| \leq B_i \}} 
-
(\Delta_i X^c)^2
\right]
+
\left[ 
\sum _{i=1}^n K_{{\delta}}(t_{i-1} - \tau) (\Delta_i X^c)^2 - \sigma_\tau^2
\right] \\
& =:
(I) + (II).
\end{split}
\end{equation}
From \eqref{eq:original_kernel_error}, we have that the second moment of (II) above converges with rate $O\left(\frac{{h}}{{\delta}}\right) + O\left({\delta}^{\varpi}\right)$. The optimal rate {of (II) is given by ${h}^{\varpi / (1 + \varpi)}$ and is attained with ${\delta} \sim {h}^{1/(\varpi + 1)}$}. Therefore, {by Proposition \ref{thm:diff_of_kw_and_tkw}, as long as $c > \varpi / (1 + \varpi)$}, (I) is of higher order than (II), in which case, (I) will be either of $o\left(\frac{{h}}{{\delta}}\right)$ or {\Blue $o\left({\delta}^{\varpi}\right)$}. This completes the proof.
\end{proof}
\begin{remark}
The leading order term {of the MSE of (\ref{eq:threshold_kernel_spot_vol})} does not depend on the {threshold}. However, by selecting the optimal threshold or its approximations, we are able to optimize the {sub-order} part of the error, which enhances the performance of the estimator in practice. Also, {since taking $c\in(2,\infty)$ does not change the asymptotic rate of convergence}, we have {\Red a} certain degree of robustness of this method.
\end{remark}

With some further assumptions, we are {also} able to obtain the {CLT} of the threshold-kernel estimator. The proof of the following {result} is similar to that Proposition \ref{prop:MSE_tkw}, {\Blue but taking advantage of Theorems 6.1 and 6.2 in \cite{FLL2016} which deal with the analogous results without jumps.}
\begin{theorem}\label{ThCLTAF}
Assume that {A}ssumption \ref{assumption:boundedness_of_intensity},
\ref{assumption:min_f_zero}, \ref{IndependentBoundednessCondition}, \ref{VolatilityConditionGeneral} and \ref{AdmissibleKernel} are satisfied, and take the threshold vector to be ${\Blue B_{n,i}(c)} = \sqrt{c \bar{\sigma}_{t_i,h}^2 {h} \log(1/{h})} + o(\sqrt{{h} \log(1 / {h})})$ for any $c\in (\frac{\varpi}{\varpi + 1}, \infty) $. Then, for each $\tau \in (0,T)$,
\begin{equation}\label{eq:CLT_Error1}
\left( \frac{{h}}{{\delta}} \right)^{-1/2} \left[ 
TKW(\tau, n, {\delta})
-
\int_0^T K_{{\delta}}(t - \tau) \sigma_t^2 dt
\right]
\rightarrow _D
\delta_1 N(0,1),
\end{equation}
{where $\delta^{2}_{1}=2\sigma_{\tau}^{4}\int K^{2}(x)dx$.} 
Furthermore, suppose that either one of the following conditions holds:
\begin{enumerate}[(1)]
\item
$\{\sigma^2_t\}_{t\geq{}0}$ is an It\^{o} process given by
$\sigma_t^2 = \sigma_{0}^{2}+\int_0^t {f_{s}} ds + \int_0^t {g_{s}} d{\Blue B_s}$, {where {\Blue $B$ is a Brownian motion independent of $W$ and} we {further assume that} $\sup_{t \in [0, T]}\mathbb{E}[|f_t|] < \infty$, $\sup_{t \in [0, T] } \mathbb{E} [ g_{t}^{2} ] < \infty$, and $\mathbb{E} [ ( g_{\tau + h} - g_\tau )^2 ] \to 0$ as $h \to 0$;}
\item 
$\sigma^2_t = f(t, Z_t)$, for a deterministic function $f:\mathbb{R} \times \mathbb{R} \to \mathbb{R}$ such that $f \in C^{1,2}(\mathbb{R})$, and a Gaussian process $\{Z_t\}_{t\geq{}0}$ satisfying Assumption \ref{VolatilityConditionGeneral} and some mild additional conditions\footnote{We refer the reader to \cite{FLL2016} for more details. In \cite{FLL2016}, we assume $\sigma_t^2 = f(Z_t)$, but it is actually trivial to generalize to the case that $\sigma_t^2 = f(t, Z_t)$ for $f \in C^{1,2}(\mathbb{R})$.}.
\end{enumerate}
Then, on an extension $(\bar{\Omega}, \bar{\mathscr{F}}, \bar{\mathbb{P}})$ of the probability space $(\Omega, \mathscr{F},\mathbb{P})$, equipped with a standard normal variable $\xi$ independent of {\Blue $\{\sigma_{t}\}_{t\geq{}0}$}, we have, for each $\tau \in (0,T)$,
\begin{equation}\label{eq:CLT_error2}
{\delta}^{-\varpi/2} \left( \int_0^T K_{{\delta}}(t - \tau) (\sigma_t^2 - \sigma_\tau^2) dt \right) \rightarrow_{D}\, \delta_2 \xi,
\end{equation}
where, under the condition (1) above, $\delta_2^2 = g(\tau, \omega)^2 \iint K(x) K(y) C(x, y) dxdy$, while, under the condition (2), $ \delta_2^2 = [f_{2}(\tau, Z_\tau)]^2 L^{(Z)}(\tau) \iint K(x)K(y) C^{(Z)}_\varpi (x,y) dxdy $. Here, $f_{2}(t, z) = \frac{\partial f}{\partial z}(t, z)$.
\end{theorem}

It is interesting to realize the {difference} between the range of $c$ allowed here and the one allowed for the integrated volatility. Indeed, for $\varpi \in (0, \infty)$, the range for spot volatility estimation is strictly larger than the range for the integrated volatility estimation. The reason is that the estimation of spot volatility is much less accurate than the integrated volatility. Therefore, we may conclude that even with a bad estimation of spot volatility, we are still able to get a threshold that is accurate enough for us to apply the threshold estimation and obtain another estimation of the spot volatility.

\subsection{Bandwidth and Kernel Selection}\label{BndwdthCal0}

With the leading order approximation we obtained from the previous subsection, we are now able to develop a feasible plug-in type bandwidth selection method. Furthermore, we can derive the optimal kernel function when the volatility is driven by Brownian motion. In this subsection, we describe all related results, which are direct {consequences} of Proposition \ref{prop:MSE_tkw}, and are parallel to results given by \cite{FLL2016}. {We refer to \cite{FLL2016} for the details of the proofs.}

The first result is the theoretical approximated optimal bandwidth, which can be obtained by taking {the derivatives} of the leading order terms in \eqref{eq:MSE_tkw} with respect to the bandwidth ${\delta}$.
\begin{proposition}\label{prop:tkw_optimal_bandwdith}
With the same assumptions as Proposition \ref{prop:MSE_tkw}, the approximated optimal bandwidth, denoted by ${\delta}^{a, opt}_n$, which is defined to minimize the leading order term of MSE in \eqref{eq:MSE_tkw}, is given by
\begin{equation}\label{eq:tkw_optimal_bandwdith}
\begin{split}
{\delta}^{a, opt}_n = &
n^{-1/(\varpi + 1)} \left[ \frac{2 T \mathbb{E}[\sigma_{\tau}^4] \int K^2(x)dx} {\varpi L(\tau) \iint K(x)K(y) C_{\varpi} (x,y) dxdy} \right] ^{1/(\varpi + 1)},
\end{split}
\end{equation}
while the {attained global minimum of the} approximated MSE is given by
\begin{equation}\label{eq:tkw_optimal_MSE}
\begin{split}
\mbox{MSE}^{a, opt}_n
& =
n^{-\varpi/(1+\varpi)} \frac{1 + \varpi}{\varpi}
\left( 2 T \mathbb{E}[\sigma_{\tau}^4] \int K^2(x)dx \right) ^{\varpi/(1+\varpi)} 
\left( \varpi L(\tau) \iint K(x)K(y) C_{\varpi} (x,y) dxdy \right) ^{1/(1+\varpi)}.
\end{split}
\end{equation}
\end{proposition}

 As shown in \cite{FLL2016}, {the resulting bandwidth obtained by replacing $ \mathbb{E}[\sigma_{\tau}^4]$ and $L(\tau)$ in the formula (\ref{eq:tkw_optimal_bandwdith}) with their integrated versions, $\int_0^T \mathbb{E}[\sigma_{\tau}^4] d\tau$ and $\int_0^T L(\tau) d\tau$,} is asymptotically equivalent to the optimal bandwidth that minimizes the integrated MSE, $\int_{0}^{T}\mathbb{E}\left[(\hat{\sigma}_{t}^{2}-\sigma_{t}^{2})^{2}\right]dt$. {\Blue In the case of a volatility process driven by Brownian motion, as in the setup (1) of Theorem \ref{ThCLTAF}, \cite{FLL2016} showed that $\varpi=1$, $C_{1}(x,y)=\min\{|x|,|y|\}{\bf 1}_{xy\geq{}0}$, and $L(t)=\mathbb{E}(g^{2}_{t})$, which leads to the formula:
 \begin{equation}\label{eq:tkw_optimal_bandwdithbb}
\begin{split}
{\delta}^{a, opt}_n = &
n^{-1/2} \left[ \frac{2 T \mathbb{E}[\int_{0}^{T}\sigma_{t}^4 dt]\int K^2(x)dx} {\mathbb{E}[\int_{0}^{T}g_{t}^{2}dt] \iint K(x)K(y) C_{1} (x,y) dxdy} \right] ^{1/2}.
\end{split}
\end{equation}
Furthermore, since, at best, we only have one realization of the path of $\sigma$ and we are working with a nonparametric setting for $\sigma$, it is natural to use  $\int_0^T \sigma_{t}^4dt$ and $\int_{0}^{T}g_{t}^{2}dt$ as  proxies of $\mathbb{E}[\int_{0}^{T}\sigma_{\tau}^4d\tau]$ and $ \mathbb{E}[\int_{0}^{T}g_{t}^{2}dt]$, respectively. These considerations suggest the following bandwidth selection method:
 \begin{equation}\label{eq:tkw_optimal_bandwdithcc}
\begin{split}
{\delta}^{a, opt}_n = &
n^{-1/2} \left[ \frac{2 T \int_{0}^{T}\sigma_{t}^4 dt\int K^2(x)dx} {\int_{0}^{T}g_{t}^{2}dt \iint K(x)K(y) C_{1}(x,y) dxdy} \right] ^{1/2}.
\end{split}
\end{equation}
Alternatively, by virtue of the independence condition in Assumption \ref{IndependentBoundednessCondition}, we can see (\ref{eq:tkw_optimal_bandwdithcc}) as an approximation of the optimal bandwidth that minimizes the conditional integrated MSE, $\mathbb{E}\left[\int_{0}^{T}(\hat{\sigma}_{t}^{2}-\sigma_{t}^{2})^{2}dt|\sigma_{s},\gamma_{s}:0\leq{}s\leq{}T\right]$. 

However, the bandwidth (\ref{eq:tkw_optimal_bandwdithcc}) is not yet {feasible}, since it depends on {the} unknown random quantities $\int_{0}^{T}\sigma_{t}^4dt$ and $\int_{0}^{T}g^{2}_{t}dt$. A well-known estimator of $\int_0^T \sigma_{t}^4 dt$ is the truncated realized quarticity, which is defined by $\widehat{IQ} =(3{h})^{-1} \sum_{i = 1}^n (\Delta_i X)^4{\bf 1}_{\{|\Delta_{i}X|<B_{n,i}\}}$ (see Proposition 1 in \cite{Mancini:2009} for consistency). The estimation of $\int_0^T g_{t}^{2}dt$ is more involved.  This quantity is sometimes called the integrated {\Red vol of vol (or vol vol for short)} and is essentially the quadratic variation of the volatility process.  \cite{FLL2016} introduced an estimator  based on the Two-time Scale Realized Quadratic Variation introduced in \cite{Two_Time_Scale}. Concretely, let} $\hat{\sigma}_{l, t_i}^2$ and $\hat{\sigma}_{r, t_i}^2$  be the left and right side estimator of $\sigma_{t_i}^2$, respectively, defined as the following:
\begin{equation}\label{eq:tkw_LeftRightKernelEstimator}
\hat{\sigma}_{l, t_i}^2
=
\frac{\sum _{j > i} K_{{\delta}}(t_{j-1} - t_i) (\Delta_j^n X)^2 \textbf{1}_{\{|\Delta_j^n X| \leq B_j\}} }
{{h} \sum _{j > i} K_{{\delta}}(t_{j-1} - \tau) \textbf{1}_{\{|\Delta_j^n X| \leq B_j\}} },
\quad
\hat{\sigma}_{r, t_i}^2
=
\frac{\sum _{j \leq i} K_{{\delta}}(t_{j-1} - t_i) (\Delta_j^n X)^2 \textbf{1}_{\{|\Delta_j^n X| \leq B_j\}} }
{{h} \sum _{j \leq i} K_{{\delta}}(t_{j-1} - \tau) \textbf{1}_{\{|\Delta_j^n X| \leq B_j\}} }.
\end{equation}
{Next,} we define the following two {\Blue finite differences:} $\Delta_i \hat{\sigma}^2 = \hat{\sigma}_{r, t_{i + 1}}^2 - \hat{\sigma}_{l, t_i}^2$, $\Delta_i^{(k)} \hat{\sigma}^2 = \hat{\sigma}_{r, t_{i + k}}^2 - \hat{\sigma}_{l, t_i}^2$. {Finally,} we can construct the following estimator:
\begin{equation}\label{eq:TSRVV}
   \widehat{IVV}_T^{\text{(tsrvv)}}= \frac{1}{k} \sum_{i = b}^{n - k - b} (\Delta_i^{(k)} \hat{\sigma}^2)^2  - \frac{n - k + 1}{nk} \sum_{i = b + k - 1}^{n - k - b} (\Delta_i \hat{\sigma}^2)^2.
\end{equation}
Here, $b$ is a small enough integer, when compared to $n$. The purpose of introducing such a number $b$ is to alleviate the boundary effect of the {one sided} estimators, since, for instance, it is expected that $\hat{\sigma}_{l, t_i}^2$ will be more inaccurate as $i$ {gets} smaller. The consistency of the TSRVV estimator can be proved by {\Blue Proposition} \ref{prop:MSE_tkw} and the corresponding results from \cite{FLL2016}.

The final result that we will mention in this subsection is about the optimal kernel function. Indeed, as was proved in \cite{FLL2016}, when the volatility is driven by Brownian motion, the optimal kernel function is given by the double exponential function.

\begin{theorem}\label{thm:tkw_optimal_kernel}
With the same assumptions as Proposition \ref{prop:MSE_tkw} and assuming $C_{\varpi}(r, s) = \min\{|r|, |s|\} \textbf{1}_{\{rs > 0\}}$, we have that the optimal kernel function that minimizes the approximated optimal MSE given by \eqref{eq:tkw_optimal_MSE} is the double exponential kernel function: 
$$K^{opt}(x) = \frac{1}{2} e^{-|x|}, \quad x \in \mathbb{R} .$$
\end{theorem}

\section{Full Implementation Scheme of The Threshold-Kernel Estimation}\label{sec:implementation}

In this section, we {propose} a complete data-driven threshold-kernel estimation scheme. {We consider several versions, depending} on whether we treat the volatility to be constant or not and whether we use {the first- or second-}order approximation formula. {One of our main interests is to investigate whether or not \emph{local} and/or \emph{second-order} thresholding {\Blue can improve} the performance of threshold estimation.}  

{Let us recall that 
the key problem at hand is jump detection; i.e., we hope to determine whether $\Delta_i N = 0$ or not. We are, of course, also interested in estimating the volatility, jump intensity, and jump density, but we are operating under the premise that effective jump detection leads to good estimation of the other model features.} 
In {Section \ref{section:optimal_thresholding}}, we {introduced} the expected number of jump misclassification as the objective function and obtained the {theoretical first and second order infill approximations} of the optimal threshold, respectively given by
\begin{equation}\label{eq:B1_B2_for_t_i}
B_{i}^{*1}
=
\left[ 
3\sigma_{i}^2 {h} \log \left( 1/ {h} \right)
\right]^{1/2}, \quad
B_{i}^{*2}
=
\sqrt{{h}} \sigma_{i} 
\left[ 
3\log \left( 1/{h} \right)
-
2 \log \left( \sqrt{2\pi} \mathcal{C}_0(f) \sigma_{i} \lambda_{i} \right)
\right]^{1/2},
\end{equation}
where, {with certain abuse of notation, we denote $\sigma_i^2 := \sigma^{2}_{t_{i}}$ and $\lambda_i := \lambda_{t_{i}}$.}
Although we have assumed that $\mathcal{C}_0(f)$ remains constant as the time evolves, we do allow non-constant volatility $\sigma_t$ and jump intensity $\lambda_t$. 

{\Blue Since estimating spot values is typically less accurate than estimating average values, {a simple first approach to implement (\ref{eq:B1_B2_for_t_i}) is to substitute $\sigma^{2}_{i}$ and $\lambda_{i}$ by their average values,  $\bar{\sigma}^{2}:=\int_{0}^{T}\sigma_{s}^{2}ds/T$ and $\bar{\lambda}:=\int_{0}^{T}\lambda_{s}ds/T$, respectively. 
This  simplification leads us to consider the following threshold sequences:}
\begin{equation}\label{eq:B1_B2_const_est}
B_{i}^{c1}
=
\left[ 
3\bar{\sigma}^{2} {h} \log \left( 1/ {h} \right)
\right]^{1/2}, \quad
B_{i}^{c2}
=
\sqrt{{h}} \bar{\sigma} 
\left[ 
3\log \left( 1/{h} \right)
-
2 \log \left( \sqrt{2\pi}\, \mathcal{C}_0(f)\, \bar{\sigma} \bar{\lambda} \right)
\right]^{1/2},
\end{equation}
where the superscript $c$  above is used to denote ``constant" volatility and jump intensity. In light of \eqref{eq:est_N_J_IV}, natural estimates of $\bar{\lambda}$ and $\bar{\sigma}^{2}$ are given by
\begin{equation}\label{eq:est_lambda_sigma}
\hat{\lambda} = 
\frac{1}{T}
\sum_{i = 1}^n \textbf{1}_{\{|\Delta_i X| > B_i\}},
\quad
\hat{\sigma}^2 =
\frac{1}{T}\sum_{i = 1}^n (\Delta_i X)^2 \textbf{1}_{\{|\Delta_i X| \leq B_i\}},
\end{equation}
respectively. 
The estimator of $\mathcal{C}_0(f)=f(0)$, as developed in Section \ref{sec:jump_density_estimation}, is given by
\begin{equation}\label{eq:threshold_kernel_est_jump_density_nonhomogeneous}
\widehat{\mathcal{C}_0(f)}:=\frac{1}{2|\{|\Delta_i X| > {B}_i\}|} \sum_{|\Delta X_i| > {B}_i} K_{{\delta}}(|\Delta_i X| - {B}_i),
\end{equation}
where the bandwidth $\delta$ is set according to Silverman's rule of thumb (\ref{eq:jump_density_bandwidth_rule_of_thumb}) and, for the threshold ${B}_i$, we could use the same threshold as in (\ref{eq:est_lambda_sigma}) or an estimate of $\widetilde{B}_{i}=\sqrt{4 {h} \sigma^2 \log(1/{h})}$ as suggested in Corollary \ref{lemma:approximate_optimal_threshold_E2}.
In the algorithms below and in the simulations of Section \ref{sec:monte_carlo_study}, we use the former threshold. Putting all together, the Algorithms \ref{algo:iterative_cnst} and \ref{algo:iterative_cnst2nd} below detail the implementation of the 1st and 2nd order constant thresholds (\ref{eq:B1_B2_const_est}).  Algorithm 1 is the same as that proposed in \cite{FLN2013} and, because it generates a nonincreasing sequence of thresholds and volatility estimates, is guaranteed to finish in finitely many steps. See the end of this section for more information about the stopping criteria for  Algorithm \ref{algo:iterative_cnst2nd}.
\begin{algorithm}
\caption{Iterative (Constant) 1st-Order Threshold Kernel Algorithm}\label{algo:iterative_cnst}
\begin{algorithmic}
\State Calculate $\hat{\sigma}^{2}_{Old}$ by \eqref{eq:est_lambda_sigma} setting $B_{i}=\infty$;
\State Initialize $B^{c1}_i=\left[ 
3\hat{\sigma}^{2} {h} \log \left( 1/ {h} \right)
\right]^{1/2}$, for $i=1,\dots, n$;
\State Calculate $\hat{\sigma}^{2}_{New}$ as in \eqref{eq:est_lambda_sigma}
with $B_{i}$ replaced with $B^{c1}_i$;
\While{{$\hat{\sigma}^{2}_{New}\neq \hat{\sigma}^{2}_{Old}$}}
\State $\hat{\sigma}^{2}_{Old}= \hat{\sigma}^{2}_{New}$;
\State Update $B^{c1}_i=\left[ 
3\hat{\sigma}^{2}_{Old} {h} \log \left( 1/ {h} \right)
\right]^{1/2}$, for $i=1,\dots, n$;
\State Calculate $\hat{\sigma}^{2}_{New}$ as in \eqref{eq:est_lambda_sigma}
with $B_{i}$ replaced with $B^{c1}_i$;
\EndWhile
\State Use final $B^{c1}_i$ for jump detection;
\end{algorithmic}
\end{algorithm}

\begin{algorithm}
\caption{Iterative (Constant) 2nd-Order Threshold Kernel Algorithm}\label{algo:iterative_cnst2nd}
\begin{algorithmic}
\State Calculate $\hat{\sigma}^{2}$ by \eqref{eq:est_lambda_sigma} setting $B_{i}=\infty$;
\State Initialize $B^{c2}_i=\left[ 
3\hat{\sigma}^{2} {h} \log \left( 1/ {h} \right)
\right]^{1/2}$, for $i=1,\dots, n$;
\While{{``Stopping Criteria" not satisfied}}
\State Calculate  $\hat{\lambda}$ and $\hat{\sigma}^{2}$ as in \eqref{eq:est_lambda_sigma}
with $B_{i}$ replaced with $B^{c2}_i$;
\State {Estimate   $\widehat{\mathcal{C}_0(f)}$ by \eqref{eq:threshold_kernel_est_jump_density_nonhomogeneous}
with $B_{i}=B_{i}^{c2}$}; 
\State Update $B^{c2}_i$ by \eqref{eq:B1_B2_const_est} with $\bar{\sigma}^{2}=\hat{\sigma}^{2}$, $\bar{\lambda}=\hat\lambda$, and $C_{0}(f)=\widehat{\mathcal{C}_0(f)}$, based on newly estimated parameters;
\EndWhile
\State Use final $B^{c2}_{i}$ for jump detection.
\end{algorithmic}
\end{algorithm}

We now consider the implementation of the local or non-constant thresholds (\ref{eq:B1_B2_for_t_i}). 
First of all, since Theorem \ref{thm:optimal_threshold_characterizations} establishes that $\sigma_i^2$ has a much greater effect on the approximated optimal threshold than that of $\lambda_i$, we simplify the problem by estimating $\lambda_{i}$ with $\hat{\lambda}$ as defined in (\ref{eq:est_lambda_sigma}). The estimation of $\sigma_i^2$, per our discussion in Section \ref{sec:threshold_kernel_spot_vol}, is given by the kernel estimator:
\begin{equation}\label{eq:threshold_kernel_spot_volbb}
\hat{\sigma}_{i}^{2}
:=
\sum _{j=1}^n K_{{\delta}}(t_{j-1} - t_{i}) (\Delta_j X)^2 \textbf{1}_{\{|\Delta_j X| \leq B_j \}}.
\end{equation}
Above, we could try to calibrate the bandwidth ${\delta}$ using an  approach similar to that described in Section \ref{BndwdthCal0}. However, for simplicity, in the simulations we set $\delta=h_{n}^{1/2}$, which, per (\ref{eq:tkw_optimal_bandwdithbb}), is rate optimal at first order.
Based on the $\hat{\sigma}_{i}^{2}$'s, $\hat{\lambda}$, and $\widehat{\mathcal{C}_0(f)}$ as defined in (\ref{eq:threshold_kernel_est_jump_density_nonhomogeneous}), we can then compute estimates of the first and second order approximation of the optimal thresholds as follows:
\begin{equation}\label{eq:B1_B2_nonconst_est}
B_{i}^{n1}
=
\left[ 
3\hat{\sigma}_i^{2} {h} \log \left( 1/ {h} \right)
\right]^{1/2}, \quad
B_{i}^{n2}
=
\sqrt{{h}} \hat{\sigma}_i 
\left[ 
3\log \left( 1/{h} \right)
-
2 \log \left( \sqrt{2\pi}\, \widehat{\mathcal{C}_0(f)}\, \hat{\sigma}_i \hat{\lambda} \right)
\right]^{1/2},
\end{equation}
where above the superscript $n$ stands for non-constant volatility estimation. Algorithm  \ref{algo:iterative_tk} below gives the details of the implementation of the non-constant thresholds (\ref{eq:B1_B2_for_t_i}). Therein, the initial threshold is taken as $B^{n1}_i=\left[ 
3\hat{\bar{\sigma}}^{2}_{0} {h} \log \left( 1/ {h} \right)
\right]^{1/2}$, where $\hat{\bar{\sigma}}^{2}_{0}$ is an initial estimate of $\bar{\sigma}^{2}:=\int_{0}^{T}\sigma_{s}^{2}ds/T$  such as those obtained from the previous Algorithms. In the simulations of Section \ref{sec:monte_carlo_study}, we take that from Algorithm \ref{algo:iterative_cnst}. See also below for more details about the ``stopping criteria" of the algorithm.
%
\begin{algorithm}
\caption{Iterative Threshold Kernel Algorithm}\label{algo:iterative_tk}
\begin{algorithmic}
\State Initialize $B^{n1}_i=\left[ 
3\hat{\bar{\sigma}}^{2}_{0} {h} \log \left( 1/ {h} \right)
\right]^{1/2}$ (or $B^{n2}_i=\left[ 
3\hat{\bar{\sigma}}^{2}_{0} {h} \log \left( 1/ {h} \right)
\right]^{1/2}$ when using 2nd order approx.), for $i=1,\dots, n$;
\While{``Stopping Criteria" not satisfied}
\State Calculate  $\hat{\sigma}_{i}^{2}$ as in \eqref{eq:threshold_kernel_spot_volbb}
with $B_{i}$ replaced with $B^{n1}_i$ (or $B^{n2}_i$ when using 2nd order approximation);
\State Calculate  $\hat{\lambda}$ and $\widehat{\mathcal{C}_0(f)}$ by \eqref{eq:est_lambda_sigma}-\eqref{eq:threshold_kernel_est_jump_density_nonhomogeneous} 
with $B_{i}$ replaced with $B^{n1}_i$ (or $B^{n2}_i$ when using 2nd order approximation);
\State Update $B^{n1}_i$ (or $B^{n2}_i$ when using 2nd order approximation) by \eqref{eq:B1_B2_nonconst_est} based on newly estimated parameters;
\EndWhile
\State Use $B^{n1}_i$ (or $B^{n2}_i$) as the final threshold value.
\end{algorithmic}
\end{algorithm}

Note} that, in \eqref{eq:B1_B2_const_est} and \eqref{eq:B1_B2_nonconst_est}, $B^{c2}$, and $B^{n2}$ may not be well defined, under a finite sample setting. {Indeed, for} a fixed time period and a fixed sample size, it is possible to have
$3\log \left( 1/{h} \right)
<
2 \log \left( \sqrt{2\pi} \mathcal{C}_0(f) \hat\sigma_{i} \hat\lambda \right)$,
in which case the square root {in (\ref{eq:B1_B2_nonconst_est})} is not well defined. Of course, asymptotically this is never an issue since we only need to consider a small enough ${h}$. 
As to implementation, {however,  it is natural to use} $B^{*2}$ whenever 
$3\log \left( 1/{h} \right)
>
2 \log \left( \sqrt{2\pi} \mathcal{C}_0(f) \sigma_{i} \lambda_{i} \right)$,
and use $B^{*1}$, otherwise.

{We now briefly discuss {some} stopping criteria for the Algorithms \ref{algo:iterative_cnst} and \ref{algo:iterative_tk}. {Typically, most iterative algorithms are stopped when} the updated value is ``close" enough to the old value. However, for the threshold estimator, we note that there are only $2^{n}$ possible threshold vectors after the initial set up. Therefore, there are only two possible situations for the ``while" loop in Algorithm \ref{algo:iterative_tk}:
\begin{enumerate}
\item
After a few iterations, the algorithm comes to a fixed  threshold vector $[\vec{B}_T]$.
\item
After a few iterations, the algorithm comes to a loop of threshold vectors given by $[\vec{B}_T^1]$, ..., $[\vec{B}_T^k]$.
\end{enumerate}
As we will see at the end of Subsection \ref{sec:MC_comparison_threshold}, generally the threshold vector converges within 2 iterations.}

\section{Monte Carlo Study}\label{sec:monte_carlo_study}
In this section, {\Blue we investigate} the performance of our proposed methods. Specifically, in Section \ref{sec:MC_comparison_threshold}, we will compare the four different threshold methods given by \eqref{eq:B1_B2_const_est} and \eqref{eq:B1_B2_nonconst_est} {\Blue and detailed in Algorithms \ref{algo:iterative_cnst}-\ref{algo:iterative_tk}}. In Section \ref{sec:MC_est_jump_density}, we investigate the performance of the threshold-kernel estimation of the jump density at the origin. 

Throughout, we consider the jump-diffusion model given by \eqref{eq:jump_diffussion_model}, {with the continuous part $\{X_t^c\}_{t \geq 0}$ following a} Heston model:
\begin{equation}\label{eq:Heston}
\begin{split}
& dX^c_t = \mu_t dt + \sqrt{V_t} dB_t , \\
& dV_t = \kappa(\theta-V_t)dt + \xi \sqrt{V_t} dW_t .
\end{split}
\end{equation}
Here, $V_t = \sigma_t^2$ is the variance process.
The parameters of \eqref{eq:Heston} are selected according to the following setting {also used in \cite{Two_Time_Scale}}: 
\begin{equation}\label{PrmAit}
	\kappa = 5,\quad  \theta = 0.04, \quad \xi = 0.5, \quad \mu_t = 0.05 - V_t / 2.
\end{equation}
As to the initial values, we use $X^c_0 = 1$ and {$V_{0}= \sigma_0^2 = 0.04$}. {\DR The unit of time in this study is 1 year and, thus, the parameter values above are annualized.} Although {\Blue  the properties} of the threshold-kernel estimators studied in this work were derived under a non-leverage setting ({\Blue i.e.,} {$\rho = 0$}, where $\rho$ is the correlation between $B_t$ and $W_t$), we run simulations on both the non-leverage setting and a negative leverage setting ($\rho = -0.5$) in order to check the robustness of the method against the leverage effect.

As to the jump component, we consider {\Blue Merton type of jumps}:
\begin{equation}\label{MertonDnsty}
{\Blue f_{normal}(x) = \frac{1}{\sqrt{2 \pi \vartheta^2}} \exp \left( -\frac{x^2}{2\vartheta^2}\right).} 
\end{equation}
The intensity of the jump component is set to be {a} constant value, i.e., $\lambda_t \equiv \lambda$ for all $t \geq 0$. For the values of $\lambda$ and $\vartheta$, we consider the following {scenarios:}
\begin{enumerate}
\item
$\lambda = 50$ {and} ${\vartheta}=0.03$, which gives an average annualized volatility of about $\sqrt{0.04+50(0.03)^{2}}\approx 0.29$;
\item
$\lambda = 100$ {and} ${\vartheta}=0.03$, which gives an average annualized volatility of about $\sqrt{0.04+100(0.03)^{2}}\approx 0.36$;

\item
$\lambda = 200$ {and} ${\vartheta}=0.03$, which gives an average annualized volatility of about $\sqrt{0.04+200(0.03)^{2}}\approx 0.46$;

\item
$\lambda = 1000$ {and} ${\vartheta}=0.01$, which gives an annualized volatility of about $\sqrt{0.04+1000(0.01)^{2}}\approx0.37$. 
\end{enumerate}
The reason for choosing these $\lambda$'s is to investigate how {\Blue large levels of} jump intensity can affect the {\Blue performance of the estimators}, {while} ${\vartheta}$ is selected accordingly {so that} the annualized volatility {\Blue is reasonable}.

{We assume} that there are 252 trading days in a year and 6.5 trading hours in each day. 
We focus on 5-minute data, which is standard in the literature to avoid microstructure noise effects. 
Furthermore, the length of the data is set to be 1 month (21 trading days), 3 month (63 trading days), and 1/2 year.

\subsection{Comparison of Different Thresholds}\label{sec:MC_comparison_threshold}

We now proceed to examine how the different ``optimal" threshold approximation methods introduced in Section 
\ref{sec:implementation} {\Blue affect} the number of jump {misclassifications}. 
{\Blue In Tables \ref{table:accuracy_misclassification_normal0}, 
 we 
report the average total number of jump mis-classifications corresponding to the four threshold approximation methods $B^{c1}$, $B^{c2}$, $B^{n1}$ and $B^{n2}$, as well as an oracle threshold, {where we use the second order approximation $B_{i}^{*2}$ in (\ref{eq:B1_B2_for_t_i})} with all the true parameter values plugged in. In each case, we} compute {\Blue the average number of jump misclassifications:}
\begin{equation}\label{eq:loss_functionSample}	
\bar{\mathcal{L}}^{a}:=\frac{1}{m}\sum_{j=1}^{m}\left(
\sum_{i=1}^{n} \textbf{1}_{ \{ |X^{(j)}_{t_{i}} - X^{(j)}_{t_{i-1}}| > B^{a}_{i,j}, N^{(j)}_{t_{i}} - N^{(j)}_{t_{i-1}} = 0 \} } 
+ 
\sum_{i=1}^{n} \textbf{1}_{ \{ |X^{(j)}_{t_{i}} - {\Red X^{(j)}_{t_{i-1}}}| \leq B^{a}_{i,j}, N^{(j)}_{t_{i}}-N^{(j)}_{t_{i-1}} \neq 0 \}} 
\right),
\end{equation}
where $m$ is the number of simulations, $X^{(j)}_{\cdot}$ and $N^{(j)}$ are the jth simulated paths of $X$ and $N$, {\Blue respectively,} and $a\in\{c1,c2,n1,n2,*2\}$, depending on the used thresholding method. {\Blue For the non-constant methods, we use an exponential kernel $K(x)=e^{-|x|}/2$ to estimate the spot volatility, which, as shown in Theorem \ref{thm:tkw_optimal_kernel}, is optimal. We ran only  4 iterations of the iterative algorithms described in Section \ref{sec:implementation}. As shown below, this typically suffices to reach convergence.}

{\Blue The conclusion is that the jump detection method based on the second order approximation with non-constant volatility estimation (``n2" method) performs the best among all the four methods. Although, {as it should be expected, this is slightly} worse than the oracle one, it is remarkably close to the latter. Even for a relatively low value of $\lambda=50$, where is typically hard to estimate $\lambda$ and $C_{0}(f)$ because of relatively few jumps, the 2nd order local method is still a bit better than those based on constant threshold. For instance, for a time horizon of 1 month, the n2 method only misses about 1 jump out of the expected 4 jumps during the month. The difference between the constant and local thresholds becomes more crucial as the intensity of jumps increases. For an intensity of $200$, the method will only miss about 3 of the expected 16 jumps.} 

{\Blue As mentioned above, the results of Table \ref{table:accuracy_misclassification_normal0} were based on 4 iterations of the Algorithms of Section \ref{sec:implementation}. To assess the convergence of the algorithm, in Table \ref{table:convergence_normal0}, we show  the results of the average number of jump misclassifications $\bar{\mathcal{L}}^{a}$, as defined in (\ref{eq:loss_functionSample}), for each of the first 4 iterations of Algorithm \ref{algo:iterative_tk} based on the 2nd order approximation. As it can be seen, convergence is typically reached after the 2nd iteration.}	

\begin {table}
\begin{center}
\begin{scriptsize}
\begin{tabular}{ | l | l | l | l | l | l | l | l | l | l | }
\hline
\#OfDays	&	\#Obs./Hr	&	$\rho$	&	$\lambda$	&	sd($f$)	&	$\bar{\mathcal{L}}^{c1}$	&	$\bar{\mathcal{L}}^{c2}$	&	$\bar{\mathcal{L}}^{n1}$	&	$\bar{\mathcal{L}}^{n2}$	&	$\bar{\mathcal{L}}^{*2}$	\\\hline
21	&	12	&	0	&	50	&	0.03	&	0.848  & 0.864	&
0.835 & 	0.795 & 0.678\\\hline
21	&	12	&	-0.5	&	50	&	0.03	&	 0.844 & 0.884	&	0.856 	&	 0.829	&	0.690	\\\hline
21	&	12	&	0	&	100	&	0.03	&	1.669	&	1.591	&	1.623	&	 1.382	&	1.259	\\\hline
21	&	12	&	-0.5	&	100	&	0.03	&	 1.628  & 1.643	&	1.584	&	 1.381	&	1.272	\\\hline
21	&	12	&	0	&	200	&	0.03	&	3.384  & 2.967	&	3.318	&	 2.603	&	2.529	\\\hline
21	&	12	&	-0.5	&	200	&	0.03	&	3.372  & 2.882 	&	3.284	&	 2.577  	&	2.487	\\\hline
21	&	12	&	0	&	1000	&	0.01	&	51.301  & 32.673 	&	49.700	&	 31.087  	&	30.218	\\\hline
21	&	12	&	-0.5	&	1000	&	0.01	&	 51.547 & 32.937 	&	49.895	&	 31.361  	&	30.480	\\\hline
63	&	12	&	0	&	50	&	0.03	&	 2.660 & 4.217 	&	2.531 	&	 2.174  	&	2.098	\\\hline
63	&	12	&	-0.5	&	50	&	0.03	&	2.590 & 4.137 	&	2.466 	&	 2.125  	&	2.051	\\\hline
63	&	12	&	0	&	100	&	0.03	& 4.952 & 6.688 & 4.741	&	3.876  	&	3.739	\\\hline
63	&	12	&	-0.5	&	100	&	0.03	& 4.914 & 6.822& 4.737	&	3.937  	&	3.820	\\\hline
63	&	12	&	0	&	200	&	0.03	&	10.195 & 11.518 & 9.842	&	7.651  	&	7.491	\\\hline
63	&	12	&	-0.5	&	200	&	0.03	&	10.001 & 11.144 & 9.658	&	7.515  	&	7.339	\\\hline
63	&	12	&	0	&	1000	&	0.01	&	148.661 & 107.325 & 143.339	&	89.477  	&	87.434	\\\hline
63	&	12	&	-0.5	&	1000	&	0.01	&	149.923 & 107.895 & 144.979	&	90.393  	&	88.293	\\\hline
126	&	12	&	-0.5	&	100	&	0.03	&	10.106 & 18.243 &9.636 & 7.890 	&	7.624	\\\hline
126	&	12	&	-0.5	&	200	&	0.03	&	20.129 & 27.433 & 19.353 & 15.036	&	14.605	\\\hline
126	&	12	&	-0.5	&	1000	&	0.01	&	298.656 & 241.770 & 285.045 & 177.588 	&	173.745	\\\hline
\end{tabular}
\end{scriptsize}
\\
\caption {{\DR Average total number of jump mis-classifications} for Normal Jumps based on 1000 Samples.}
\label{table:accuracy_misclassification_normal0}
\end{center}
\end {table}


\begin {table}
\begin{center}
\begin{scriptsize}
\begin{tabular}{ | l | l | l | l | l | l | l | l | l  | }
\hline
\#OfDays	&	\#Obs./Hr	&	$\rho$	&	$\lambda$	&	sd($f$)	&	$\bar{\mathcal{L}}^{n2,{\rm Iter 1}}$	&	$\bar{\mathcal{L}}^{n2,{\rm Iter 2}}$	&	$\bar{\mathcal{L}}^{n2,{\rm Iter 3}}$	&	$\bar{\mathcal{L}}^{n2,{\rm Iter 4}}$		\\\hline
21	&	12	&	0	&	50	&	0.03	&	0.795 & 0.795
 & 0.795 & 0.795  \\\hline
21	&	12	&	-0.5	&	50	&	0.03	&	 0.830  & 0.829 &0.829 & 0.829\\\hline
21	&	12	&	0	&	100	&	0.03	&	1.383 & 1.382 & 1.382 & 1.382	\\\hline
21	&	12	&	-0.5	&	100	&	0.03	&	 1.384 & 1.382  & 1.381 & 1.381		\\\hline
21	&	12	&	0	&	200	&	0.03	&	2.602 & 2.603 & 2.603
& 2.603		\\\hline
21	&	12	&	-0.5	&	200	&	0.03	&	2.586 & 2.577 & 2.577
 & 2.577 		\\\hline
21	&	12	&	0	&	1000	&	0.01	&	31.855 & 31.216
&  31.121 & 31.087   	\\\hline
21	&	12	&	-0.5	&	1000	&	0.01	&	 32.080 & 31.482
& 31.385 & 31.361  		\\\hline
63	&	12	&	0	&	50	&	0.03	&	 2.183 & 2.174  & 2.174&  2.174  		\\\hline
63	&	12	&	-0.5	&	50	&	0.03	&	2.126 & 2.125 & 2.125 & 2.125 		\\\hline
63	&	12	&	0	&	100	&	0.03	& 3.875 & 3.874 & 3.876 &  3.876 		\\\hline
63	&	12	&	-0.5	&	100	&	0.03	& 3.952 & 3.938 & 3.937
 & 3.937	\\\hline
63	&	12	&	0	&	200	&	0.03	&	7.680 & 7.653 & 7.651 & 7.651		\\\hline
63	&	12	&	-0.5	&	200	&	0.03	&	7.544 & 7.517 & 7.515
 &  7.515 		\\\hline
63	&	12	&	0	&	1000	&	0.01	&	91.283
&   89.747  & 89.526  & 89.477  		\\\hline
63	&	12	&	-0.5	&	1000	&	0.01	&	92.324
&   90.722  & 90.411  & 90.393  		\\\hline
126	&	12	&	-0.5	&	100	&	0.03	&	7.929  & 7.892 
&  7.889  & 7.890\\
\hline
126	&	12	&	-0.5	&	200	&	0.03	&15.097 & 15.036 &
  15.035 & 15.036	 	\\\hline
126	&	12	&	-0.5	&	1000	&	0.01	& 181.403
 & 178.022 & 177.677 & 177.588	 	\\\hline
\end{tabular}
\end{scriptsize}
\\
\caption {Average total number of jump mis-classifications for the first 4 iterations of the nonhomogeneous Algorithm \ref{algo:iterative_tk} based on 2nd order approximations. Again, we use 1000 Samples for all values of $T$.}
\label{table:convergence_normal0}
\end{center}
\end {table}

\subsection{Estimation of Jump Density at the Origin {\Blue and Spot Volatility}}\label{sec:MC_est_jump_density}

{We now study} the performance of the kernel estimator of the jump density at the origin that we proposed in Section \ref{sec:jump_density_estimation}. Since we have already confirmed that the second order approximation of the optimal threshold with non-constant volatility estimation outperforms other thresholds, we will only consider {this threshold in this} and later subsections.

{\Blue The results are shown in {\Blue Table \ref{table:accuracy_normal_density_uniformaa}}. These basically confirm what we expect that the performance of the estimator improves as the time-horizon and intensity become larger (for the same level of jump variance). It is hard to compare the performance of the estimators when $\vartheta=0.03$ to those when $\vartheta=0.01$ and $\lambda=1000$ because, though we expect more jumps in the latter case, those will also be much harder to detect since $\vartheta$ is smaller.} Finally, an interesting phenomenon}  is that we usually underestimate the jump density at the origin. This is acceptable for our purpose. Indeed, if we denote $\widehat{B^2}$ as the estimated second order threshold, we generally have $B^1 > \widehat{B^2} > B^2$. This is better than having $\widehat{B^2} < B^2$, in which case we might suffer significantly from {false positives (i.e., mis-classifying the increments of the continuous component as jumps)}.

\begin {table}
\begin{center}
\begin{scriptsize}
\begin{tabular}{ | l | l | l | l | l | l | l | l | l | }
\hline
\#OfDays	&	\#Obs./Hr	&	$\rho$	&	$\lambda$	&	$\vartheta={\rm sd}(f)$	&	$f(0)$	&	${\rm E}(\hat{f}^{n2}(0))$	&	${\rm sd}(\hat{f}^{n2}(0))$	&	{\Blue $\sqrt{MSE(\hat{f}^{n2}(0))}$}		\\\hline
21	&	12	&	0	&	100	&	0.03	&	13.30	&	$ 8.6965$	&	$5.4797$	&	$7.1567$	\\\hline
21	&	12	&	-0.5	&	100	&	0.03	&	13.30	&	$8.8117$	&	$5.6950$	&	$7.2510$	\\\hline
63	&	12	&	0	&	100	&	0.03	&	13.30	&	$11.6005$	&	$2.1668$	&	$2.7537$	\\\hline
63	&	12	&	-0.5	&	100	&	0.03	&	13.30	&	$11.4145$	&	$2.2175$	&	$2.9107$	\\\hline
126	&	12	&	-0.5	&	100	&	0.03	&	13.30	&	$11.9558$	&	$1.6286$	&	$2.1116$	\\\hline
21	&	12	&	0	&	200	&	0.03	&	13.30	&	$11.2759$	&	$2.6362$	&	$3.3236$	\\\hline
21	&	12	&	-0.5	&	200	&	0.03	&	13.30	&	$11.1900$	&	$2.7161$	&	$3.4393$	\\\hline
63	&	12	&	0	&	200	&	0.03	&	13.30	&	$11.9714$	&	$1.6997$	&	$2.1582$	\\\hline
63	&	12	&	-0.5	&	200	&	0.03	&	13.30	&	$11.9234$	&	$1.6539$	&	$2.1518$	\\\hline
126	&	12	&	-0.5	&	200	&	0.03	&	13.30	&	$12.4776$	&	$1.3081$	&	$1.5451$	\\\hline
21	&	12	&	0	&	1000	&	0.01	&	39.89	&	$37.9363$	&	$4.5286$	&	$4.9321$	\\\hline
21	&	12	&	-0.5	&	1000	&	0.01	&	39.89	&	$37.8582$	&	$4.2041$	&	$4.6693$	\\\hline
63	&	12	&	0	&	1000	&	0.01	&	39.89	&	$41.4335$	&	$2.9176$	&	$3.3007$	\\\hline
63	&	12	&	-0.5	&	1000	&	0.01	&	39.89	&	$41.5071$	&	$3.0726$	&	$3.4722$	\\\hline
126	&	12	&	-0.5	&	1000	&	0.01	&	39.89	&	$41.8874$	&	$2.4255$	&	$3.1420$	\\\hline
\end{tabular}
\end{scriptsize}
\\
\caption {MSE of Jump Density Estimation at the origin $0$ for Normal Jumps based on 1000 Samples.}
\label{table:accuracy_normal_density_uniformaa}
\end{center}
\end {table}

{\Blue We finally give some illustrations about the performance of the the kernel/threshold spot volatility estimator (\ref{eq:threshold_kernel_spot_volbb}). We apply 4 iterations of the local Algorithm \ref{algo:iterative_tk} based on the 2nd order approximation of the optimal threshold. In Figure \ref{ApproxNum1}, we show a prototypical realization of the variance process $\{V_{t}\}_{t\geq{}0}$ defined in (\ref{eq:Heston}) together with the estimated spot variance process resulting from the 1st iteration (red dotted), from the  final iteration 4 (long-dashed blue), and from the oracle (green double-dashed), which uses $B_{i}^{*1}$ in (\ref{eq:B1_B2_for_t_i}) with the true values of $\sigma_{i}^{2}=V_{t_{i}}$, $C_{0}(f)$, and $\lambda$. We take $\lambda=200$, $\vartheta=0.03$, $T=6$ months, and $h=5$ minutes. The three spot variance estimates are close to each other and are able to fit well the overall level of the volatility through time. The Sum Of Square Errors,
\[
	{\rm SSE}=\sum_{i=1}^{n}(\hat{\sigma}_{t_{i}}^{2}-\sigma^{2}_{t_{i}})^{2},
\]
for the 1st, 4th, and oracle estimates are respectively given by $1.6525$, $1.4457$, and $1.4450$. Figure \ref{ApproxNum2} shows the same results corresponding to $\lambda=1000$ and $\vartheta=0.01$. The SSE are in this case $2.0015$, $1.5069$, and $1.4006$ for the 1st, 4th, and oracle estimates, respectively.}

\begin{figure}[h]
\begin{center}
	\includegraphics[width=15.0cm,height=10cm]{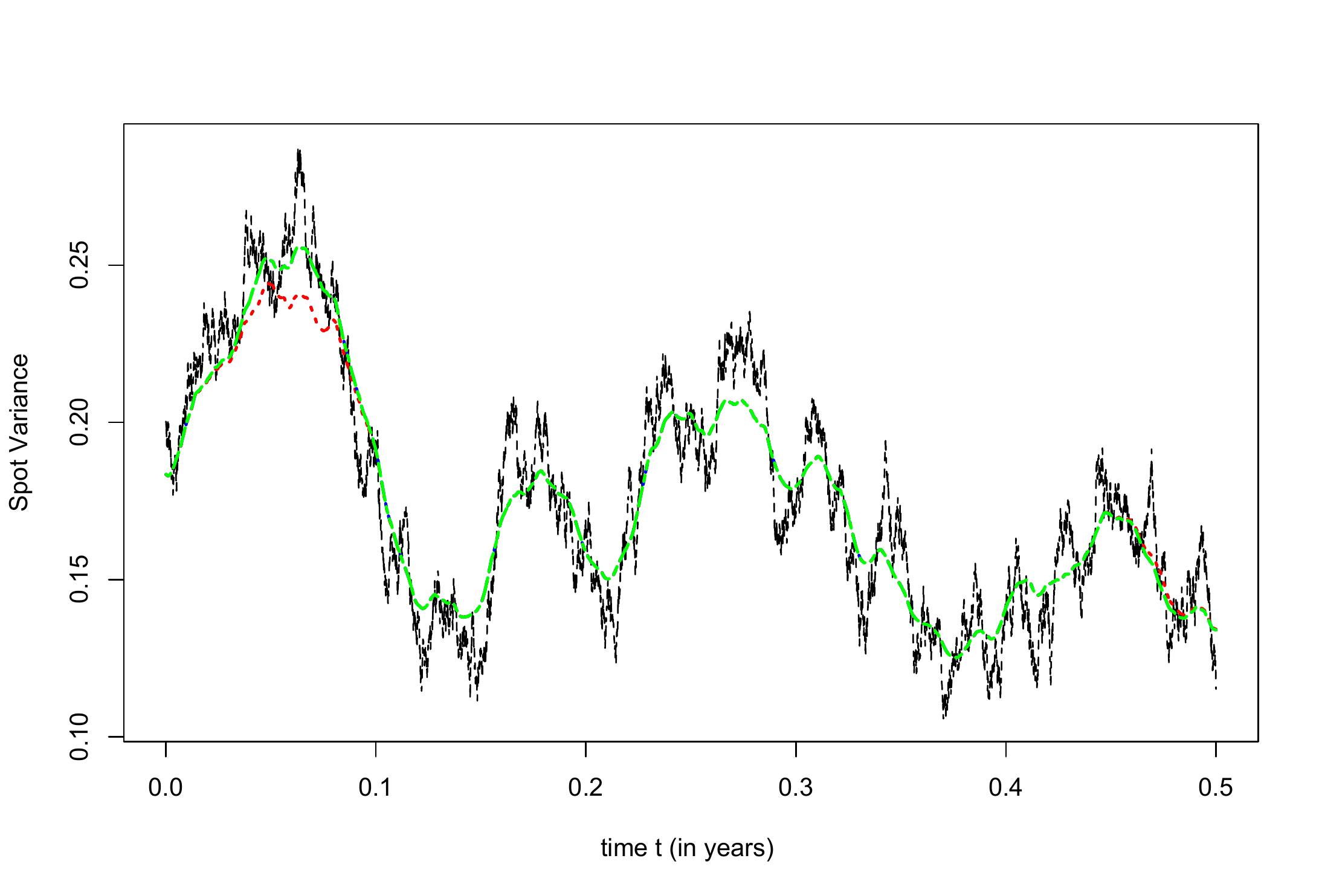}
\end{center}
    \caption{Variance process $\{V_{t}\}_{t\geq{}0}$ (jiggling dotted black curve) and the estimated spot variance process (\ref{eq:threshold_kernel_spot_volbb}) resulting from the 1st iteration (red dotted) and 4th iteration 4 (long-dashed blue) of Algorithm \ref{algo:iterative_tk} based on $B^{n2}$. We also plot the oracle variance process (\ref{eq:threshold_kernel_spot_volbb}) (green double-dashed) replacing $B_{i}$ with the true $B_{i}^{*}$ in (\ref{eq:B1_B2_for_t_i}). The oracle and the 4th iteration variance estimates overlap. We take the parameter values in (\ref{PrmAit}) as well as $\rho=-0.5$, $\lambda=200$, and the Merton jumps (\ref{MertonDnsty}) with $\vartheta=0.03$. Estimation based on 5-minute observations during 6 months.}\label{ApproxNum1}
\end{figure}

\begin{figure}[h]
\begin{center}
	\includegraphics[width=15.0cm,height=10cm]{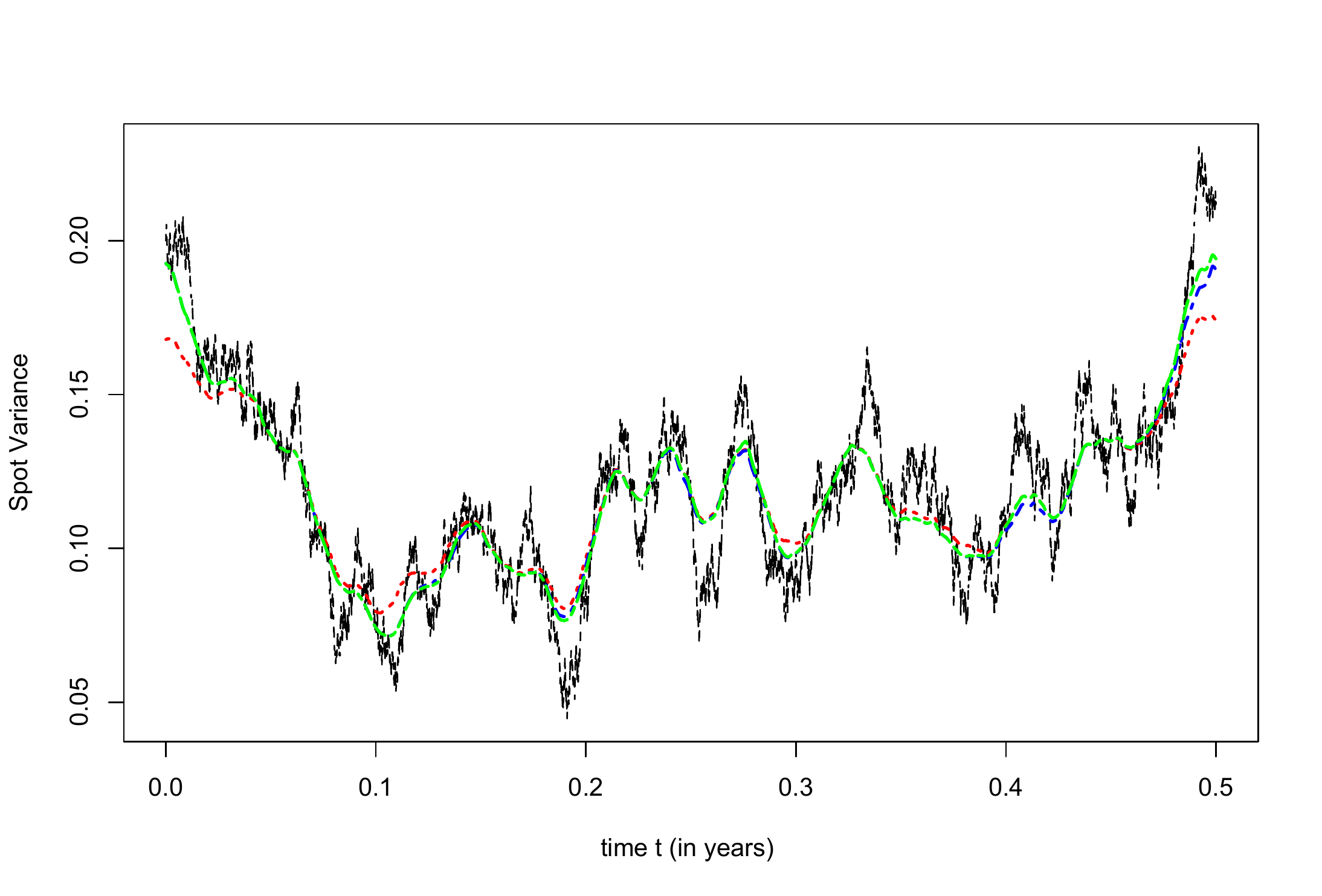}
\end{center}
    \caption{Variance process $\{V_{t}\}_{t\geq{}0}$ (jiggling dotted black curve) and the estimated spot variance process (\ref{eq:threshold_kernel_spot_volbb}) resulting from the 1st iteration (red dotted) and 4th iteration 4 (long-dashed blue) of Algorithm \ref{algo:iterative_tk} based on $B^{n2}$. We also plot the the oracle variance process (\ref{eq:threshold_kernel_spot_volbb}) (green double-dashed) replacing $B_{i}$ with the true $B_{i}^{*}$ in (\ref{eq:B1_B2_for_t_i}). The oracle and the 4th iteration variance estimates almost overlap during the whole domain except at the end. We take the parameter values in (\ref{PrmAit}) as well as $\rho=-0.5$, $\lambda=1000$, and the Merton jumps (\ref{MertonDnsty}) with $\vartheta=0.01$. Estimation based on 5-minute observations during 6 months.}\label{ApproxNum2}
\end{figure}

\section{Conclusion and Future Work}\label{sec:conclusion}

In this paper, we study the problem of jump detection via the thresholding method, which {is obviously closely} related to the problem of spot volatility estimation. We extend the approximated optimal threshold of \cite{FLN2013} by considering a second-order approximation and a non-homogeneous parameter setting. The result is of theoretical interest since the remainder of the second order approximation is much smaller and, at the same time, the resulting threshold estimator is time-invariant, which makes more sense in reality. Monte Carlo studies also demonstrate the superior performance of the second-order approximation.

The higher accuracy comes with the price of more parameters to estimate. We first managed to build a threshold-kernel estimator of the jump density at the origin. We propose a different ``optimal" threshold for this purpose and demonstrate the reason why this should be different from the original ``optimal" threshold. The intuition is that we have to be more accurate when claiming that an increment contains a jump in order to have a good estimation of its density at the origin. We also put forward a modified version of the threshold-kernel estimator of spot volatility where increments that exceed the threshold are filtered out.

In order to implement the proposed methods, we need to resolve some key {obstacles}. Concretely, {estimates of} the optimal threshold, the jump density at the origin, and the spot volatility depend on each other. To resolve the issue, we propose an iterative threshold-kernel estimation scheme. Although we are not guaranteed that the iterative algorithm always converges, Monte Carlo studies show that this rarely creates any problem in reality.

{The spirit of jump detection by threshold method is to claim that a jump occurs whenever the absolute value of the increment of the process exceeds the threshold, which, by definition, is a binary outcome. In this case, when an increment is {close} to the threshold, a small difference in the increment can lead to totally different results. One way to alleviate such a problem is to estimate the probability that a jump happens during a specific time interval, which is similar to the idea of Logistic regression. This {\Red suggests} an alternative approach to threshold-based classification. 
Given a non-decreasing function $F:[0, \infty) \to [0, 1]$ and an increment $|\Delta_{i}X|$, we can postulate that the probability that a jump occurs during $[t_{i-1},t_{i}]$ is $F(|\Delta_{i} X|)$. We can then adopt the following loss function, that is frequently {used in} classification problems:
\begin{align*}
L_{t,h}(F) 
= 
\Ex \left( F(|X_{t+h}-X_{t}|) \mathbf{1}_{\{N_{t+h}-N_{t} = 0\}} \right) 
 + 
\Ex \left( [ 1 - F(|X_{t+h}-X_{t}|)] \mathbf{1}_{\{N_{t+h}-N_{t} \neq 0\}} \right).
\end{align*}
Indeed, it could be cumbersome  to optimize over all continuous functions $F$. However, we can try to limit ourself to a suitable, relatively small, class of possible functions $F$. One possible direction is to consider $F_B(x) = F(x  / B)$, which is a generalization of what we have done in this paper. Another possible direction is to consider $F(x) = F_{1n}(x) \mathbf{1}_{\{ 0 < x < B\}} + F_{2n}(x) \mathbf{1}_{\{ x \geq B\}}$, where $F_{1n}$ and $F_{2n}$ are two functions that can depend on $n$. This can potentially provide insight on how the shape of $F$ should look like around the ``optimal" threshold.}

\appendix

\section{Proof of the Main Results}  \label{sec:appendix_optimal_thresholding}

Let us start by giving a lemma necessary  for the proof of Theorem \ref{thm:uniform_convexity}.
\begin{lemma} \label{lemma:sufficient_condition2} For $i = 1, 2$ let $f_{i} \in \mathcal{C}\left([0,\infty)\right)$  be strictly positive and differentiable on $(0,\infty)$.  Further suppose that $f_{1}$ is non-increasing while $f_{2}$ is non-decreasing and $lim_{x \to{} 0^{+}} [f_{1}^{'}(x) + f_{2}^{'}(x)]$ exists.  If there exists $x_{0} \in (0,\infty)$ such that 
\begin{equation}
(a) \;\; |f^{'}_{1}(x)| \; \geq \; |f^{'}_{2}(x)| \quad \text{for all} \;\; x \in (0, x_{0}) \hspace{10mm}
(b) \;\; |f^{'}_{2}(x)| \; \geq \; |f^{'}_{1}(x)| \quad \text{for all} \;\; x \in (x_{0}, \infty), \label{eq:quasi_convex_diff_condition}
\end{equation}
then, $f := f_{1} + f_{2}$ is quasi-convex on $[0,\infty)$.
\end{lemma}
\begin{proof}
From $(a)$ in (\ref{eq:quasi_convex_diff_condition}) and since $f_{1}^{'}(x) \leq 0$ and $f_{2}^{'}(x) \geq 0$, $f^{'}(x) = f_{1}^{'}(x) + f_{2}^{'}(x) \leq 0$ for all $x \in (0,x_{0})$.  Furthermore, this implies $\lim_{x \to{} 0^{+}} f^{'}(x) \leq 0$.  On the other hand, from $(b)$ in  (\ref{eq:quasi_convex_diff_condition}) $f^{'}(x) = f_{1}^{'}(x) + f_{2}^{'}(x) \geq 0$ for all $x \in (x_{0},\infty)$.  From the well known sufficient conditions for the quasi-convexity of continuous real-valued functions of a real variable ({see} \cite{BV2004} pg. 99 (3.20) therein for more details), it follows that $f$ is quasi-convex on $[0,\infty)$.
\end{proof}

\begin{proof} [\textbf{Proof of Theorem \ref{thm:uniform_convexity}}]
{\Blue Throughout, we assume $w=1$ and $\overline{\gamma}_{t,h} > 0$ (the cases of $\overline{\gamma}_{t,h} \leq 0$ and $w\neq{}0$} can be proved in a similar way). {\Blue For simplicity we omit the argument $w$ in $L_{t,h}(B;w)$}. Let $F^{*k}_{t,h}$ denotes the distribution of the density $\phi_{t,h}*f^{*k}$.
Conditioning on the number of jumps $N_{t+h}-N_{t}$, the loss function $L_{t,h}$ is split as follows:
\[
L_{t,h}(B) :=  L_{t,h}^{(1)}(B) + L_{t,h}^{(2)}(B),
\]
where
\begin{align*}
	L_{t,h}^{(1)}(B) 
	&:= 
	\Px \left( |X_{t+h}-X_{t}| > B , N_{t+h}-N_{t} = 0 \right)
	=
	e^{- h\overline{\lambda}_{t,h}} \left[ 
	1 - \Phi \left( \frac{B - h\overline{\gamma}_{t,h}}{\overline{\sigma}_{t,h} \sqrt{h}} \right) 
	+
	\Phi \left( \frac{-B - h \overline{\gamma}_{t,h}}{\overline{\sigma}_{t,h} \sqrt{h}} \right) 
	\right], \\
	L_{t,h}^{(2)}(B) 
	&:= 
	\Px \left( |X_{t+h}-X_{t}| \leq B , N_{t+h}-N_{t} \neq 0 \right)
	=
	e^{- h \overline{\lambda}_{t,h}} \sum_{k=1}^{\infty} \frac{\left( h \overline{\lambda}_{t,h}\right)^{k}}{k!} \left[ F_{t,h}^{*k}(B) - F_{t,h}^{*k}(-B) \right].
\end{align*}
Here, $\Phi(\cdot)$ is the cdf of standard normal distribution. Note that by definition, $L_{t,h}^{(1)}$ is strictly decreasing while $L_{t,h}^{(2)}$ is strictly increasing.  It is also clear that for each $h > 0$ and $t \in [0,T]$, $L_{t,h}^{(1)} \in \mathcal{C}^{\infty}\left(\mathbb{R}^{+}\right)$ and $\partial _B L_{t,h}^{(1)} (B) < 0$ for all $B \in \mathbb{R}^+$. For the differentiability of $L_{t,h}^{(2)}$,  since 
\begin{equation} \label{eq:uniform_density_bound}
\sup_	{k \in \mathbb{N}} \; \displaystyle{\sup_{x \in \mathbb{R}}} \left| \phi_{t,h}*f^{*k}(x) \right| \leq \sup_{x\in\Rx}f(x)=:M(f) < \infty,
\end{equation}
it follows that $ \sup_{k \in \mathbb{N}}\sup_{B \in (0,\infty)} \left| \phi_{t,h}*f^{*k}(B) + \phi_{t,h}*f^{*k}(-B) \right| \leq 2 M(f)$ and, thus, by Bounded Convergence Theorem, $L_{t,h}^{(2)}$ is differentiable. {Similarly, since $ \sup_{m \in \mathbb{N}} \sup_{k \in \mathbb{N}}\sup_{B \in (0,\infty)} \left| \phi^{(m)}_{t,h}*f^{*k}(B) + \phi^{(m)}_{t,h}*f^{*k}(-B) \right| \leq 2 M(f)$, we can further prove that $L_{t,h}^{(2)} \in \mathcal{C}^{\infty}\left(\mathbb{R}^{+}\right)$ by Bounded Convergence Theorem.}
%
%
%
%
%
%

{We observe that $L_{t,h}^{(1)}(B) \neq 0$ and $L_{t,h}^{(2)}(B) \neq 0$ for all $B > 0$, so}
we now proceed to study the ratio
\[
	R_{t,h}(B) := \frac{ \partial_{B}L_{t,h}^{(2)}(B)}{- \partial_{B}L_{t,h}^{(1)}(B)}.
\]
Let us start by noting that
\begin{align}
	\partial_{B}L_{t,h}^{(1)}(B) &= - \frac{e^{-h\overline{\lambda}_{t,h}}}{ \sqrt{h}\overline{\sigma}_{t,h}} \left[ {\phi} \left( \frac{B - h\overline{\gamma}_{t,h}}{\overline{\sigma}_{t,h} \sqrt{h}} \right) + {\phi} \left( \frac{B + h\overline{\gamma}_{t,h}}{\overline{\sigma}_{t,h} \sqrt{h}} \right) \right],
\label{eq:loss_deriv_comp1}
	\\
	\partial_{B}L_{t,h}^{(2)}(B)&= e^{- h\overline{\lambda}_{t,h}} \sum_{k=1}^{\infty} \frac{\left( h\overline{\lambda}_{t,h}\right)^{k}}{k!} \left[ \phi_{t,h}*f^{*k}(B) + \phi_{t,h}*f^{*k}(-B) \right].
\label{eq:loss_derive_comp2}
\end{align}
An immediate consequence is that $R_{t, h}(B)$ is continuous for $B \in[0, \infty)$. $R_{t, h}$ may now be written as:
\[
R_{t,h}(B) = \sum_{k=1}^{\infty} \frac{\left( h\overline{\lambda}_{t,h}\right)^{k}}{k!} I_{t,h,k}(B),\quad\text{where}
\quad I_{t,h,k}(B) 
:= 
\frac{ \overline\sigma_{t,h} \sqrt{h}\left(\phi_{t,h}*f^{*k}(B)+\phi_{t,h}*f^{*k}(-B)\right)}{\phi \left( \frac{B-h\overline{\gamma}_{t,h}}{\overline{\sigma}_{t,h}\sqrt{h}} \right) + \phi \left( \frac{B+h\overline{\gamma}_{t,h}}{\overline{\sigma}_{t,h}\sqrt{h}} \right)}.
\]
By definition of convolution, $I_{t,h,k}$ can be written as: 
\begin{align*}
I_{t,h,k}(B) 
&=
\int g_{t,h}(w,B) f^{*k}(w)dw,\quad \text{where}\quad g_{t,h}(w,B) 
:= 
\frac{\phi \left( \frac{B- h\overline{\gamma}_{t,h}-w}{\overline{\sigma}_{t,h} \sqrt{h}} \right)
+ 
\phi \left( \frac{B+ h\overline{\gamma}_{t,h}+w}{\overline{\sigma}_{t,h} \sqrt{h}} \right)
}{\phi \left( \frac{B- h\overline{\gamma}_{t,h}}{\overline{\sigma}_{t,h} \sqrt{h}} \right)
+
\phi \left( \frac{B+ h\overline{\gamma}_{t,h}}{\overline{\sigma}_{t,h} \sqrt{h}} \right)}.
\end{align*}

Plugging in the normal p.d.f., $g_{t,h}$ can be factorized to be:
\begin{align*}
g_{t,h}(w,B) 
= 
\exp \left(-\frac{w^{2} + 2w h\overline{\gamma}_{t,h}}{2 h\overline{\sigma}^{2}_{t,h}} \right)
\frac{\exp \left( B( h\overline{\gamma}_{t,h} + w)/ h\overline{\sigma}^{2}_{t,h} \right)
+ 
\exp \left( -B( h\overline{\gamma}_{t,h} + w)/ h\overline{\sigma}^{2}_{t,h} \right) } 
{ \exp \left( B \overline{\gamma}_{t,h}/\overline{\sigma}^{2}_{t,h} \right)
+ 
\exp \left( -B \overline{\gamma}_{t,h}/\overline{\sigma}^{2}_{t,h} \right) }
=:g_{t,h}^{(1)}(w) \; g_{t,h}^{(2)}(w,B).
\end{align*}
It is not hard to prove the following properties of $g_{t,h}^{(1)}$:
\begin{equation}\label{eq:Thm21_proof_g1_bound}
1 \leq g_{t,h}^{(1)}(w) \leq e^{ h\overline{\gamma}^{2}_{t,h}/2\overline{\sigma}^{2}_{t,h}},\quad 
\omega \in (-2 h\overline{\gamma}_{t,h},0), 
\quad \mbox{ and }\quad
0 < g_{t,h}^{(1)}(w) \leq 1,\quad 
\omega \in (-2 h\overline{\gamma}_{t,h},0)^C.
\end{equation}
%
%
$g^{(2)}_{t, h}(\omega, \cdot)$ is a function of type 
$t(x) = \frac{e^{ax} + e^{-ax}}{e^{bx} + e^{-bx}}$
, where $x\in [0,\infty)$, $a = |( h\overline{\gamma}_{t,h} + w)/ h\overline{\sigma}^{2}_{t,h}|$ and $b = |\overline{\gamma}_{t,h}/\overline{\sigma}^{2}_{t,h}|$.
Note that the derivative $t^\prime (x)$ can be written as
$$
t^\prime(x) 
=
\frac{e^{ax} + e^{-ax}}{e^{bx} + e^{-bx}}
\left[
a\frac{e^{ax} - e^{-ax}}{e^{ax} + e^{-ax}}
-
b\frac{e^{bx} - e^{-bx}}{e^{bx} + e^{-bx}}
\right].
$$
%
When $a > b > 0$, $t(x)$ is an increasing function from $1$ to $+\infty$ and 
\begin{equation*}
\begin{split}
t^\prime(x)
& \geq
\frac{e^{ax}}{2e^{bx}}
(a - b)
\frac{e^{bx} - e^{-bx}}{e^{bx} + e^{-bx}}
\geq
\frac{(a - b)}{4}e^{(a - b)x}
( 1 - e^{-2bx} )
\geq
\left\{
\begin{array}{ll}
\frac{a - b}{4}(1 - e^{-1})2bx, & x\leq \frac{1}{2b} \\
\frac{a - b}{4}(1-e^{-1})(a-b)x , & x > \frac{1}{2b} \\
\end{array} 
\right. \\
& \geq
\frac{(a - b)(1 - e^{-1})\min(a - b, 2b)}{4}x.
\end{split}
\end{equation*}
{For the third inequality, when $0 \leq x \leq 1/2b$, we use $1 - e^{-2bx} \geq (1 - e^{-1})2bx$, and when $x > 1/2b$, we use $e^{(a - b)x} \geq (a - b)x$ and $1 - e^{-2bx} \geq (1 - e^{-1})$.}
Specifically, when $a > 3b$, we have 
\begin{equation}\label{eq:Thm21_proof_g2_a_geq_b}
t^\prime(x) \geq b^2(1 - e^{-1})x.
\end{equation}
When $b > a > 0$, $t(x)$ is a decreasing function from $1$ to $0$ and 
\begin{equation}\label{eq:Thm21_proof_g2_a_less_b}
\left| t^\prime(x) \right|
\leq
b\frac{e^{bx} - e^{-bx}}{e^{bx} + e^{-bx}}
\leq b^2x,
\end{equation}
{where we use the property that $\tanh^\prime(x) \leq 1$.}

Here we notice that $a < b \Leftrightarrow \omega \in (-2 h \overline{\gamma}_{t,h},0)$.
Based on this, for each fixed $k \in \mathbb{N}$, we decompose $I_{t,h,k}$ into two parts:
\begin{equation}\label{eq:Thm21_proof_decompose_I}
I_{t,h,k}(B) 
= 
\left( \int_{(-2 h \overline{\gamma}_{t,h},0)} + \int_{(-2 h \overline{\gamma}_{t,h},0)^{c}} \right)
g_{t,h}^{(1)}(w) \; g_{t,h}^{(2)}(w,B) f^{*k}(w)dw
=: 
I_{t,h,k}^{(1)}(B) + I_{t,h,k}^{(2)}(B).
\end{equation}

In what follows, We shall prove that there exists $h_0 > 0$, which may depend on $T$, such that for all $t \in [0, T]$ and $h \in (0, h_0)$, there exists $B^{*}_{t,h}>0$, such that 
\begin{equation*}
	R_{t,h}(B) < 1,\quad \text{for }B\in(0,B^{*}_{t,h}), \quad\text{ and }\quad
	 R_{n}(B) > 1,\quad \text{for }B\in(B^{*}_{t,h},\infty).
\end{equation*}
These two conditions, together with the signs of $\partial_{B}L_{t,h}^{(1)}$ and $\partial_{B}L_{t,h}^{(2)}$, will imply that $B\to{}L_{t,h}(B)$ is 
quasi-convex (see Lemma \ref{lemma:sufficient_condition2} below) 
for $h$ small enough.
To do this, we will prove the following:
\begin{enumerate}[(i)]
\item
For any $h > 0$, $\lim_{B \to \infty} R_{t,h}(B) = +\infty$.
\item
$\lim_{h \to 0} \sup_{t \in [0, T]} R_{t,h}(0) = 0$.
\item
There exists $h_{0} > 0$, which may depend on $T$, such that for all $t \in [0, T]$ and $h \in (0, h_0)$, $R_{t,h}(\cdot)$ is strictly increasing.
\end{enumerate}
  
For (1), it is clear that $I_{t,h,k}^{(1)} \geq 0$, and by Fatou's Lemma, for $k$ large enough
\footnote{{$k$ has to be large, since now we are not assuming small $h$, so it is possible that $f^{*k}(\omega) \equiv 0$ for $\omega \in (-2 h \overline{\gamma}_{t,h},0)^{c}$.
}}
, $I_{t,h,k}^{(2)}$ satisfies
$$
\liminf_{B \to \infty} I^{(2)}_{t,h,k}(B) 
\geq
\int_{(-2 h \overline{\gamma}_{t,h},0)^{c}} 
\liminf_{B \to \infty} g_{t,h}^{(1)}(w) \; g_{t,h}^{(2)}(w,B) f^{*k}(w)dw
=
+\infty.
$$
{These two relationships  imply (i)}.

{For (ii)}, since $g_{t,h}^{(2)}(w,0) = 1$,
\begin{align*}
I_{t,h,k}(0) 
= 
\int g^{(1)}_{t,h}(w) f^{*k}(w)dw 
= 
\sqrt{2\pi h}\overline{\sigma}_{t,h} e^{h\overline{\gamma}_{t,h}^2/2\overline{\sigma}_{t,h}^2}
\int \frac{ e^{-(w + h\overline{\gamma}_{t,h})^{2}/2 h\overline{\sigma}_{t,h}^2}}{\sqrt{2\pi h} \overline{\sigma}_{t,h}} f^{*k}(w)dw 
\leq 
\sqrt{2\pi h}\overline{\sigma}_{t,h} e^{h\overline{\gamma}_{t,h}^2/2\overline{\sigma}_{t,h}^2} M(f).
\end{align*}
Note that the right-hand side converges to zero as $h \to 0$, and does not depend on $k$. By Assumption {\ref{assumption:boundedness_of_intensity}}, the convergence is uniformly in $t$, so (ii) follows.

Now we proceed to consider (iii). Indeed, for any given $t \in [0, T]$, by the upper bound of {\Blue $g_{t,h}^{(1)}( w)$} given by \eqref{eq:Thm21_proof_g1_bound} and the upper bound of $\left| \partial_B g_{t,h}^{(2)}( w,B) \right|$ given by \eqref{eq:Thm21_proof_g2_a_less_b}, we have
\begin{align*}
|I_{t,h,k}^{(1)}(B+\delta) - I_{t,h,k}^{(1)}(B)|
& =
\int_{(-2 h \overline{\gamma}_{t,h},0)}
g_{t,h}^{(1)}(w) \times | g_{t,h}^{(2)}(w,B + \delta) - g_{t,h}^{(2)}(w,B) | \times f^{*k}(w)dw \\
& \leq
\int_{(-2 h \overline{\gamma}_{t,h},0)}
e^{ h\overline{\gamma}^{2}_{t,h}/2\overline{\sigma}^{2}_{t,h}} 
\times 
\frac{\overline{\gamma}_{t,h}^2}{\overline{\sigma}^{4}_{t,h}} (B + \delta) \delta 
\times f^{*k}(w)dw
\leq
2h
\frac{\overline{\gamma}_{t,h}^3}{\overline{\sigma}^{4}_{t,h}} 
e^{ h\overline{\gamma}^{2}_{t,h}/2\overline{\sigma}^{2}_{t,h}} M(f) 
(B + \delta) \delta .
\end{align*}
Furthermore, for $I_{t,h,k}^{(2)}$, note that for $\omega \in (-2 h \overline{\gamma}_{t,h}, 0)^C$, $g_{t,h}^{(2)}(w,B)$ is increasing in $B$, and for $\omega \in [-4 h \overline{\gamma}_{t,h}, 4 h \overline{\gamma}_{t,h}]^C$, we have $|( h\overline{\gamma}_{t,h} + w)/ h\overline{\sigma}^{2}_{t,h}| > 3 |\overline{\gamma}_{t,h}/\overline{\sigma}^{2}_{t,h}|$. Thus, we have
\begin{align*}
I_{t,h,k}^{(2)}(B+\delta) - I_{t,h,k}^{(2)}(B)
& =
\int_{(-2 h \overline{\gamma}_{t,h},0)^C}
g_{t,h}^{(1)}(w) \times ( g_{t,h}^{(2)}(w,B + \delta) - g_{t,h}^{(2)}(w,B) ) \times f^{*k}(w)dw \\
& \geq
\int_{[-4 h \overline{\gamma}_{t,h}, 4 h \overline{\gamma}_{t,h}]^C}
g_{t,h}^{(1)}(w) 
\times 
\frac{(1 - e^{-1})\overline{\gamma}^2_{t,h}}{\overline{\sigma}^{4}_{t,h}}
B\delta 
\times 
f^{*k}(w)dw \\
& \geq
\frac{(1 - e^{-1})\overline{\gamma}^2_{t,h}}{\overline{\sigma}^{4}_{t,h}}
B\delta 
\left( \int
g_{t,h}^{(1)}(w) f^{*k}(w)dw - 8h\overline{\gamma}_{t,h} e^{ h\overline{\gamma}^{2}_{t,h}/2\overline{\sigma}^{2}_{t,h}} M(f) \right).
\end{align*}
Putting these two inequalities together, we have that for any $B > 0$ and $0 < \delta < B$:
\begin{align*}
\frac{1}{B\delta} \sum_{k=1}^{\infty} \left|
 \frac{\left( h\overline{\lambda}_{t,h}\right)^{k}}{k!} \left( I_{t,h,k}^{(1)}(B+\delta) - I_{t,h,k}^{(1)}(B) \right)
 \right|
& \leq
4h \left( e^{h\overline{\lambda}_{t,h}} - 1 \right) 
\frac{\overline{\gamma}_{t,h}^3}{\overline{\sigma}^{4}_{t,h}}
e^{ h\overline{\gamma}^{2}_{t,h}/2\overline{\sigma}^{2}_{t,h}} M(f) = O(h^2), 
\quad h \rightarrow 0,
\end{align*}
and
\begin{align*}
& \frac{1}{B\delta}
\sum_{k=1}^{\infty}
\frac{\left( h\overline{\lambda}_{t,h}\right)^{k}}{k!} 
\left( I_{t,h,k}^{(2)}(B+\delta) - I_{t,h,k}^{(2)}(B) \right) \\
&\quad \geq
\frac{(1 - e^{-1})\overline{\gamma}^2_{t,h}}{\overline{\sigma}^{4}_{t,h}}
\left( h\overline{\lambda}_{t,h} \int
g_{t,h}^{(1)}(w) f(w)dw - 
\sum_{k=1}^{\infty}
\frac{\left( h\overline{\lambda}_{t,h}\right)^{k}}{k!} 8h\overline{\gamma}_{t,h} e^{ h\overline{\gamma}^{2}_{t,h}/2\overline{\sigma}^{2}_{t,h}} M(f) \right) \\
&\quad \geq
h^{3/2} \overline{\lambda}_{t,h}
\frac{(1 - e^{-1})\overline{\gamma}^2_{t,h}}{\overline{\sigma}^{4}_{t,h}}
\sqrt{2\pi}\overline{\sigma}_{t,h}
\exp \left( \frac{h\overline{\gamma}^2_{t,h}}{2 \overline{\sigma}^{2}_{t,h}} \right) \frac{\mathcal{C}_m(f)}{2} + O(h^2),
\quad h \rightarrow 0,
\end{align*}
where the last equality can be justified by $\int g_{t,h}^{(1)}(w) f(w)dw \geq \sqrt{2\pi h}\overline{\sigma}_{t,h}
\exp \left( \frac{h\overline{\gamma}^2_{t,h}}{2 \overline{\sigma}^{2}_{t,h}} \right) \frac{\mathcal{C}_m(f)}{2} + O(h)$ for small $h$, where $\mathcal{C}_m(f)$ is defined in \eqref{eq:density_notation}, since the following holds: 
$$
g_{t,h}^{(1)}(w)
=
\exp \left(-\frac{w^{2} + 2w h\overline{\gamma}_{t,h}}{2 h\overline{\sigma}^{2}_{t,h}} \right)
=
\sqrt{2\pi h}\overline{\sigma}_{t,h}
\exp \left( \frac{h\overline{\gamma}^2_{t,h}}{2 \overline{\sigma}^{2}_{t,h}} \right)
\frac{1}{\sqrt{2\pi h}\overline{\sigma}_{t,h}}
\exp \left(- \frac{(w + h\overline{\gamma}_{t,h})^2 }{2 h\overline{\sigma}^{2}_{t,h}} \right).
$$ 
{Also note} that both convergences do not depend on $B$ and $\delta$, and by Assumption {\ref{assumption:boundedness_of_intensity}}, both the convergences can all be made uniform in $t$. This proves (iii).
\end{proof}

\begin{proof} [\textbf{Proof of Theorem \ref{thm:optimal_threshold_characterizations}}]
For simplicity, we use the notation $f_{t,h}^{*k}:=\phi_{t,h}*f^{*k}$, {\Blue where recall that $
\phi_{t,h}(x):= \frac{1}{\overline{\sigma}_{t,h}\sqrt{h}}\phi \left( \frac{x- h\overline{\gamma}_{t,h}}{\overline{\sigma}_{t,h}\sqrt{h}} \right)
$ is  the density of $X^{c}_{t+h}-X^{c}_{t}$.} We start by demonstrating that the optimal thresholds $(B^{*}_{t,h})_{t,h}$ converge to $0$ uniformly on $t \in [0,T]$, as $h \to{} 0$. {Let us first note that the loss function {(\ref{eq:loss_single_term})} can be written as
\begin{align*}
	L_{t,h}(B) &:= e^{-h\overline\lambda_{t,h}}\Px \left( \left|h\overline\gamma_{t,h}+\overline\sigma_{t,h}\sqrt{h}Z\right|> B \right) +  e^{-h\bar\lambda_{t,h}}\sum_{k=1}^{\infty} \frac{(h\overline\lambda_{t,h})^{k}}{k!}\Px \left( \Big|h\overline\gamma_{t,h}+\overline\sigma_{t,h}\sqrt{h}Z+\sum_{i=1}^{k}\zeta_{i}\Big| \leq B \right).
\end{align*}
{\Blue Next, by partitioning $E:=\{|h\overline\gamma_{t,h}+\overline\sigma_{t,h}\sqrt{h}Z+\sum_{i=1}^{k}\zeta_{i}| \leq B \}$ into $E\cap \{|h\overline\gamma_{t,h}+\overline\sigma_{t,h}\sqrt{h}Z|\leq B\}$ and 
$E\cap \{|h\overline\gamma_{t,h}+\overline\sigma_{t,h}\sqrt{h}Z|> B\}$ and simplifying,} 
\begin{align*}
		L_{t,h}(B) &\leq{}\Px \left( \left|h\overline\gamma_{t,h}+\overline\sigma_{t,h}\sqrt{h}Z\right|> B \right) +e^{-h\bar\lambda_{t,h}}\sum_{k=1}^{\infty} \frac{(h\overline\lambda_{t,h})^{k}}{k!}\Px \left( 
		{\Blue E}, \Big|h\overline\gamma_{t,h}+\overline\sigma_{t,h}\sqrt{h}Z\Big|\leq B  \right)\\
				&\leq{}\Px \left(h\gamma^{*}_{T+h}+\sigma^{*}_{T+h}\sqrt{h}|Z|> B \right) +\sum_{k=1}^{\infty} \frac{(h\lambda^{*}_{T+h})^{k}}{k!}\Px \left( \Big|\sum_{i=1}^{k}\zeta_{i}\Big| \leq 2B \right),
\end{align*}
where we have used that $\gamma_{t}^{*}:=\sup_{s\leq{}t}|\gamma_{s}|$, $\sigma_{t}^{*}:=\sup_{s\leq{}t}\sigma_{s}$, and $\lambda_{t}^{*}:=\sup_{s\leq{}t}\lambda_{s}$ are finite for any $t$.
{\Blue Next, consider} a sequence of thresholds given by $B_{h,c}^{Pow} := ch^{\alpha}$ for $\alpha \in (0,1/2)$ and $c > 0$. Thus, using that $\Px \left( \left|\zeta_{1}\right| \leq 2B \right)\sim 4\mathcal{C}_{0}(f)B$ and $\Px \left(|\sum_{i=1}^{k}\zeta_{i}| \leq 2B \right)=O(B)$ as $B\to{}0$, 
\[
	\sup_{t\in[0,T]}L_{t,h}\left(B^{Pow}_{h,c}\right)
	\leq 
	4\,c\,\mathcal{C}_{0}(f)\lambda^{*}_{T+h}h^{1+\alpha}+o(h^{1+\alpha}).
\]
Now suppose that $\epsilon := \limsup_{h \to{} 0^{+}} \sup_{t \in [0,T]} B^{*}_{t,h} > 0$.  Then, there exists subsequences $(h_{n})_{n}$ and $(t_{n})_{n}$ such that $\inf_{n} B_{t_{n},h_{n}}^{*} \geq \epsilon/2$. In that case, 
\[
	L_{t_{n},h_{n}}(B_{t_{n,h_{n}}}^{*})\geq{} e^{-h{\lambda}^{*}_{T+h}}h\underline\lambda_{T+h}\Px \left( \left|h_{n}\overline\gamma_{t_{n},h_{n}}+\overline\sigma_{t_{n},h_{n}}\sqrt{h_{n}}Z+\zeta_{1}\right| \leq\epsilon/2\right),
\]
but, also $L_{t_{n},h_{n}}(B_{t_{n,h_{n}}}^{*})\leq L_{t_{n},h_{n}}(B_{h_{n},c}^{Pow})$ and, since $\Px \left( \left|h_{n}\overline\gamma_{t_{n},h_{n}}+\overline\sigma_{t_{n},h_{n}}\sqrt{h_{n}}Z+\zeta_{1}\right| \leq\epsilon/2\right) {\to \Px \left( \left| \zeta_{1}\right| \leq\epsilon/2\right) > 0}$, 
\footnote{It is necessary to have $\mathcal{C}_0 (f) > 0$. Otherwise, $B \to 0$ is not optimal.}
as $n\to\infty$, we would have that 
\[
	 4\,c\,\lambda^{*}_{T+h}\mathcal{C}_{0}(f)h^{1+\alpha}+o(h^{1+\alpha})\geq h\underline\lambda_{T+h}+o(h),
\]
which leads to a contradiction.}
Hence, it is necessary that the optimal thresholds converge to $0$ uniformly on $[0,T]$. 

Now we will show the asymptotic characterization of the optimal thresholds.  From Theorem \ref{thm:uniform_convexity}, there exists $h_{0} > 0$, {depending on $T$,} such that, for all $t \in [0,T]$ and $h \in (0,h_{0}]$, the loss functions {$L_{t,h}$} possess a unique critical point.  
By equating the first-order derivative of the loss function to zero, from (\ref{eq:loss_deriv_comp1})-(\ref{eq:loss_derive_comp2}) it follows that the unique optimal threshold, $B_{t,h}^{*}$, must satisfy the equation given by
\begin{equation}  \label{eq:grad_zero_Ito_semi}
\frac{1}{{\sqrt{h}\overline{\sigma}_{t,h}}} \left[ {\phi} \left( \frac{B_{t,h}^{*} - {h\overline{\gamma}_{t,h}}}{{\sqrt{h}}\overline{\sigma}_{t,h}} \right) + {\phi} \left( \frac{B_{t,h}^{*} + {h}\overline{\gamma}_{t,h}}{{\sqrt{h}}\overline{\sigma}_{t,h}} \right) \right] = \sum_{k=1}^{\infty} \frac{\left({h}\overline{\lambda}_{t,h}\right)^{k}}{k!} \left[ f_{t,h}^{*k}(B_{t,h}^{*}) + f_{t,h}^{*k}(-B_{t,h}^{*}) \right].
\end{equation}
A rearrangement of this equation shows
\begin{equation}  \label{eq:opt_thresh_zero}
{\phi} \left( \frac{B_{t,h}^{*} - {h}\overline{\gamma}_{t,h}}{{\sqrt{h}}\overline{\sigma}_{t,h}} \right) = {\sqrt{h}}\overline{\sigma}_{t,h} \left[ 1 + e^{-2B_{t,h}^{*}\overline{\gamma}_{t,h}/\overline{\sigma}^{2}_{t,h}} \right]^{-1}  \sum_{k=1}^{\infty} \frac{\left({h}\overline{\lambda}_{t,h}\right)^{k}}{k!} \left[ f_{t,h}^{*k}(B_{t,h}^{*}) + f_{t,h}^{*k}(-B_{t,h}^{*}) \right].
\end{equation}
Upon taking the log on both sides of (\ref{eq:opt_thresh_zero}), we arrive at the fixed point equation (\ref{eq:fixed_point_equation}).
From (\ref{eq:uniform_density_bound}) together {\Blue with Assumption \ref{assumption:boundedness_of_intensity}}, we conclude that $\lim_{h \to{} 0^{+}} B_{t,h}^{*}/h^{1/2} = \infty$, uniformly for $t \in [0,T]$, {i.e.,} 
\[
	\lim_{h\to{}0}\inf_{t\in[0,T]} \frac{B_{t,h}^{*}}{\sqrt{h}}=+\infty.
\]  
A further modification of this equation indicates that
\begin{align*}
B_{t,h}^{*} &= {h}\overline{\gamma}_{t,h} + \sqrt{2{h}} \overline{\sigma}_{t,h} \log^{1/2} \left( \frac{1}{\overline{\sigma}_{t,h}\overline{\lambda}_{t,h}{h^{3/2}}}\right)
\left[ 
1 
+ 
\frac{\log \left(\frac{\sqrt{2\pi} \left(f_{t,h}^{*1}(B_{t,h}^{*})+f_{t,h}^{*1}(-B_{t,h}^{*})\right)}{1 + e^{-2B_{t,h}^{*}\overline{\gamma}_{t,h}/\overline{\sigma}^{2}_{t,h}}}\right)}{\log \left(\overline{\sigma}_{t,h} \overline{\lambda}_{t,h}{h^{3/2}} \right)} 
+ 
\frac{\log \left( 1 + S_{t,h}(B_{t,h}^{*}) \right)}{\log \left(\overline{\sigma}_{t,h} \overline{\lambda}_{t,h} {h^{3/2}}\right)} 
\right]^{1/2},
\end{align*}
where above, we have defined
\[
 S_{t,h}(B) :=  \sum_{k=2}^{\infty} \frac{\left({h}\overline{\lambda}_{t,h}\right)^{k-1}}{k!} \left[ \frac{f_{t,h}^{*k}(B) + f_{t,h}^{*k}(-B)}{f_{t,h}^{*1}(B) + f_{t,h}^{*1}(-B)} \right].
\]
From this, a direct consequence is that $B^{*}_{t, h} = O ( \sqrt{h\log(1/h)})$, so 
we have
\begin{equation*}
\frac{\sqrt{2\pi} \left(f_{t,h}^{*1}(B_{t,h}^{*})+f_{t,h}^{*1}(-B_{t,h}^{*})\right)}{1 + e^{-2B_{t,h}^{*}\overline{\gamma}_{t,h}/\overline{\sigma}^{2}_{t,h}}}
=
\sqrt{2\pi} \mathcal{C}_0 (f) + O(B^*_{t, h}), \quad
S_{t,h}(B_{t,h}^{*}) = O(h^2).
\end{equation*}
The second relationship above is because $f_{t, h}^{*k}$ are bounded by $M(f)$ and {\Blue $f_{t, h}^{*1}(B_{t,h}^{*})$} is bounded away from zero. We prove the first relationship above now. {\Blue Indeed, 
by} our assumption on the smoothness of $f$, there exists $\epsilon > 0$, such that $f \in C^1((0, \epsilon))$ and $f \in C^1((-\epsilon , 0))$. Then, we have:
\begin{equation}\label{eq:Thm22_proof_convolution_difference}
\begin{split}
f_{t,h}^{*1} (B^*_{t, h}) &= f_{t,h}^{*1}(0) + O(B^*_{t, h})\\
f_{t,h}^{*1}(0) - \mathcal{C}_0(f)
& =
f_{t,h}^{*1}(0) - 
\left( 
f(0-) \int_{-\infty}^0 \phi_{t, h}(y)  dy
+
f(0+) \int_0^{+\infty} \phi_{t, h}(y) dy 
\right) 
+ O(\sqrt{h}) \\
& =
\int_{-\infty}^0 (f(y) - f(0-)) \phi_{t, h}(y)  dy +\int_0^{+\infty} (f(y) - f(0+))\phi_{t, h}(y) dy  + O(\sqrt{h}) \\
& =
\int_{-\epsilon}^0 (f(y) - f(0-)) \phi_{t, h}(y)  dy +\int_0^{+\epsilon} (f(y) - f(0+))\phi_{t, h}(y) dy  + O(\sqrt{h}) \\
& =
\int_{-\epsilon}^0{\Blue \int_{0}^{1} f^\prime(yv) dvy} \phi_{t, h}(y)  dy +\int_0^{+\epsilon}{\Blue  \int_{0}^{1}f^\prime(yv)dv} y\phi_{t, h}(y)  dy  + O(\sqrt{h}) \\
& =
O(\sqrt{h}).
\end{split}
\end{equation}
Above, the first equality uses $\int_0^{+\infty} \phi_{t, h}(y) dy = 1/2 + O(\sqrt{h})$ and $\int_{-\infty}^0 \phi_{t, h}(y) dy = 1/2 + O(\sqrt{h})$. The third equality uses $\int_\epsilon^{+\infty} \phi_{t, h}(y) dy = o(h)$ and {\Blue $\int_{-\infty}^{-\epsilon} \phi_{t, h}(y) dy = o(h)$. 
From this, we have $f_{t,h}^{*1}(0) = \mathcal{C}_0(f) + O(\sqrt{h})$. We then have $f_{t,h}^{*1}(B_{t,h}^{*})+f_{t,h}^{*1}(-B_{t,h}^{*}) = 2 \mathcal{C}_0(f) + O(B_{t,h}^{*})$.
Therefore,} for any $\alpha \in (0, 1/2)$,
\begin{align*}
\frac{\log \left(\frac{\sqrt{2\pi} \left(f_{t,h}^{*1}(B_{t,h}^{*})+f_{t,h}^{*1}(-B_{t,h}^{*})\right)}{1 + e^{-2B_{t,h}^{*}\overline{\gamma}_{t,h}/\overline{\sigma}^{2}_{t,h}}}\right)}{\log \left(
\overline{\sigma}_{t,h} \overline{\lambda}_{t,h} h^{3/2}
\right)} 
=
\frac{\log \left( \sqrt{2\pi} \mathcal{C}_0 (f) \right) }
{\log \left( \overline{\sigma}_{t,h} \overline{\lambda}_{t,h} h^{3/2} \right)} 
+ o(h^{\alpha}), \quad
\frac{\log \left( 1 + S_{t,h}(B_{t,h}^{*}) \right)}{\log \left(\overline{\sigma}_{t,h} \overline{\lambda}_{t,h} h^{3/2}\right)} 
=
o(h^{\alpha}){\Blue .}
\end{align*}
{\Blue For the last assertion of the theorem, if we further note that 
$\overline{\sigma}_{t,h}^{2} = \sigma_{t}^{2} + O(h)$ and
$\overline{\lambda}_{t,h} = \lambda_{t} + O(h)$ under the specified smoothness of $t\to\sigma_{t}^2$ and $t\to\lambda_{t}$, then} we conclude the following approximation of $B^*_{t, h}$:
\begin{align*}
B_{t,h}^{*} 
& = 
\sqrt{2h} \overline{\sigma}_{t,h} \log^{1/2} \left( \frac{1}{ \overline{\sigma}_{t,h} \overline{\lambda}_{t,h} h^{3/2} }\right)
\left[ 
1 
+ 
\frac{\log \left( \sqrt{2\pi} \mathcal{C}_0 (f) \right) }{\log \left( \overline{\sigma}_{t,h} \overline{\lambda}_{t,h} h^{3/2} \right)}
\right]^{1/2} + o(h^{\frac{1}{2} + \alpha}) \\
& =
\sqrt{h} {\Blue \sigma_{t}} 
\left[ 
3\log \left( 1/h \right)
-
2 \log \left( \sqrt{2\pi} \mathcal{C}_0(f) {\Blue \sigma_{t} \lambda_{t}} \right)
\right]^{1/2}
+ 
o(h^{\frac{1}{2} + \alpha} )
,
\end{align*}
for any $\alpha \in (0, 1/2)$. 
\end{proof}

\begin{proof}
[\textbf{Proof of Proposition \ref{prop:conditions_for_convergence_E2}}]
First, note that 
\begin{equation}\label{eq:Prop31_proof_DeltaX_geq_B}
\begin{split}
P(|\Delta X| > B)
& =
e^{-{h} \lambda} \Px \left( |{h} \gamma + \sqrt{{h}} \sigma Z | > B \right) 
+
{h} \lambda e^{-{h} \lambda} \Px \left( \left| {h} \gamma +  \sqrt{{h}}\sigma Z + \zeta \right| > B \right)
+
O({h}^2).
\end{split}
\end{equation}
{Let $\phi_{{h}}(x)$ be
the density of ${h} \gamma +  \sqrt{{h}}\sigma Z$ and note that, {for $k\geq{}1$,} ${h} \gamma +  \sqrt{{h}}\sigma Z + \sum_{i = 1}^k \zeta_{i}$ has density $\phi_{{h}} * f^{*k}$, which is bounded by $M(f):=\sup_{x}f(x)$.
Therefore, we have
$
\Big|\frac{\partial}{\partial B}
\Px \left( | {h} \gamma +  \sqrt{{h}}\sigma Z + \sum _{i = 1}^k \zeta _i | > B \right)\Big| < 2M(f)
$
and, furthermore,}
\begin{equation*}
\frac{\partial}{\partial B}
\sum_{k \geq 2} \frac{( {h} \lambda )^k}{k!} \Px \left(\Big| {h} \gamma +  \sqrt{{h}}\sigma Z + \sum _{i = 1}^k \zeta _i \Big| > B \right)
=
\sum_{k \geq 2} \frac{( {h} \lambda )^k}{k!} 
\frac{\partial}{\partial B} 
\Px \left( \Big| {h} \gamma +  \sqrt{{h}}\sigma Z + \sum _{i = 1}^k \zeta _i \Big| > B \right)
=
O({h}^2).
\end{equation*}
{We then have}
\begin{equation}\label{eq:Prop31_proof_fstar}
f^*(B) = {- \frac{\partial}{\partial B} P( |\Delta X| > B)}=
e^{-{h} \lambda} [ \phi_{{h}}(B) + \phi_{{h}}(-B)]
+
{h} \lambda e^{-{h} \lambda} [ g(B) + g(-B)]
+ O({h} ^2)
,
\end{equation}
where $g$ denotes the density of ${h} \gamma +  \sqrt{{h}}\sigma Z + \zeta$. Combining \eqref{eq:Prop31_proof_DeltaX_geq_B} and \eqref{eq:Prop31_proof_fstar}{,} the conditional density is {such that}
\begin{align*}
f^*_{ |\Delta X| | |\Delta X| > B}(B)
=
\frac{f^*(B)}{P( |\Delta X| > B)}
=
\frac{
\frac{1}{\lambda} \frac{1}{\sqrt{2 \pi {h}^3 \sigma^2}} \left[ \exp\left(  -\frac{(B - {h} \gamma)^2}{2 {h} \sigma^2 }\right) + \exp\left(  -\frac{( B + {h} \gamma)^2}{2 {h} \sigma^2 }\right) \right]
+
g(B) + g(-B)
}{
\frac{1}{\lambda {h}} \Px \left( | {h} \gamma + \sqrt{{h}} \sigma Z | > B \right) 
+
\Px \left( \left| {h} \gamma +  \sqrt{{h}}\sigma Z + \zeta \right| > B \right)
}
+ O({h}).
\end{align*}
Now, by $g = \phi_{h} * f$ and {the} smoothness of $f$ near $0$, we have that if $x$ is close enough to $0$,
\begin{align*}
g(x) - f(x) 
& = 
\left(\int_{(x - \epsilon, x + \epsilon)} +  \int_{(x - \epsilon, x + \epsilon)^c} \right)
(f(y) - f(x)) \phi_{h}(y - x) dy \\
& = 
\int_{(x - \epsilon, x + \epsilon)} (f^\prime(x) (y - x) + f^{\prime\prime} (\theta_y) (y - x)^2) \phi_{h}(y - x) dy 
+ o({h})
=
O({h}),
\end{align*}
where $\theta_y$ is between $x$ and $y$ and $\epsilon$ is a fixed positive number such that $f \in C^2((x - \epsilon, x + \epsilon))$. Such an $\epsilon$ exists due to Assumption \ref{assumption:smooth_density}. Above, we have used the following facts: 
\begin{align*}
& \int_{(x - \epsilon, x + \epsilon)^C} \phi_{h}(y - x) dy = o({h}), \quad
\int_{(x - \epsilon, x + \epsilon)} (y - x) \phi_{h}(y - x) dy = \gamma {h} + o({h}), \\
& {\int_{(x - \epsilon, x + \epsilon)} \left|f^{\prime\prime} (\theta_y)\right| (y - x)^2 \phi_{h}(y - x) dy} \leq M \sigma^2 {h}.
\end{align*}
Note that the above holds uniformly in $x$ near $0$, so we have $g(B) = f(B) + O({h})$ for ${h}$ small enough. 
This also implies
$$
\Px \left( \left| {h} \gamma +  \sqrt{{h}}\sigma Z + \zeta \right| \leq B \right)
=
2g(0)B + o(B)
=
2f(0)B + o(B).
$$
Therefore, we have the following:
\begin{align*}
f^*_{ |\Delta X | | | \Delta X | > B}(B)
& =
\frac{
\frac{1}{\lambda} \frac{1}{\sqrt{2 \pi {h}^3 \sigma^2}} \left[ \exp\left(  -\frac{(B - {h} \gamma)^2}{2 {h} \sigma^2 }\right) + \exp\left(  -\frac{( B + {h} \gamma)^2}{2 {h} \sigma^2 }\right) \right]
+
2f(0) + O({h}) + O(B^2)
}{
\frac{1}{\lambda {h}} \Px \left( | {h} \gamma + \sqrt{{h}} \sigma Z | > B \right) 
+
1 - 2f(0)B + o(B)
}
+ O({h}) \\
& =
2f(0)
+
\frac{ 2 }{\lambda \sqrt{2 \pi {h}^3 \sigma^2}} \exp\left(  -\frac{B^2}{2 {h} \sigma^2 }\right)
+
2f(0) B +  o(B)+o(h^{-3/2}e^{-\frac{B^{2}}{2h\sigma^{2}}}),
\end{align*}
where we used the following:
\begin{align*}
& \frac{1}{{h}} \Px \left( | {h} \gamma + \sqrt{{h}} \sigma Z | > B \right) 
\sim 
\frac{\sigma}{B\sqrt{2 \pi {h}}} \exp\left(  -\frac{B^2}{2 {h} \sigma^2 } \right), \quad
\frac{1}{\sqrt{2 \pi {h}^3 \sigma^2}} \exp\left(  -\frac{(B - {h} \gamma)^2}{2 {h} \sigma^2 }\right) 
\sim
\frac{1}{\sqrt{2 \pi {h}^3 \sigma^2}} \exp\left(  -\frac{B^2}{2 {h} \sigma^2 } \right).
\end{align*}
This completes the proof.
\end{proof}

\begin{proof}[\textbf{Proof of Corollary \ref{lemma:approximate_optimal_threshold_E2}}]
Denote the leading order term of \eqref{eq:kernel_jump_density_E2} as:
$$
F(B) =
\frac{ 1 }{\lambda  f(0) \sqrt{2 \pi {h}^3 \sigma^2}} \exp\left(  -\frac{B^2}{2 {h} \sigma^2 } \right)
+ B.
$$
Set $a = 1/( \lambda  f(0) \sqrt{2 \pi {h}^3 \sigma^2})$, $b = 1/(2 {h} \sigma^2)$. For ${h}$ small enough, we do have $a\sqrt{b} > 1/(1 - \exp(-1/2)) $, and $\log(2ab) < b$. By the Lemma \ref{DmyLma} below, the minimum of $F$ is in $\left( \sqrt{2{h} \sigma^2}, \sqrt{2{h} \sigma^2 \log(1/\sqrt{2 \pi {h}^5 \sigma^6})} \right)$ and satisfies $B\exp\left(  -\frac{B^2}{2 {h} \sigma^2 } \right) = \sqrt{2 \pi {h}^5 \sigma^6}$.
Taking {logarithms} on both sides and rearranging terms, we get
$$
\frac{B^2}{2{h} \sigma^2} = \log (B) - \frac{5}{2} \log ({h}) + C,
$$ 
for some constant $C$.
Note that {since $B$ lies in $\left( \sqrt{2{h} \sigma^2}, \sqrt{2{h} \sigma^2 \log(1/\sqrt{2 \pi {h}^5 \sigma^6})} \right)$, $\log (B) = \frac{1}{2} \log({h}) + O(\log\log(1/{h}))$. Thus,} we get the approximation of the optimal $B$ as
\begin{equation*}
B^* = \sqrt{4 {h} \sigma^2 \log(1/{h})} + O(\sqrt{{h} \log\log(1/{h})}).
\end{equation*}
This completes the proof.
\end{proof}

\begin{lemma}\label{DmyLma}
Suppose $a, b > 0$ and $a\sqrt{b} > 1/(1 - \exp(-1/2)) $, and $\log(2ab) < b$. Define $F(x) = a\exp(-bx^2) + x$ where $x\geq 0$. Then, the minimum point of $F$ is in $(1/\sqrt{2b}, \sqrt{\log(2ab)/b})$ and satisfies $2abx\exp(-bx^2) = 1$.
\end{lemma}
\begin{proof}
Taking derivative twice, we get $F^\prime(x) = -2abx\exp(-bx^2) + 1$ and $F^{\prime\prime}(x) = 2ab(2bx^2 - 1) \exp(-bx^2)$.
By studying the sign of $F^{\prime\prime}$, we have that $F^{\prime}$ is decreasing in $(0, 1/\sqrt{2b})$ and increasing in $(1/\sqrt{2b}, \infty)$, and we also have $F^{\prime}(1/\sqrt{2b}) = -a\sqrt{2b}\exp(-1/2) + 1$.
Now since $a\sqrt{2b} > 1/(1 - \exp(-1/2))  > \exp(1/2)$, $F^\prime(0) = F^\prime(+\infty) = 1$, we have that $F^\prime$ has a root $r_1$ in $(0, 1/\sqrt{2b})$ and another root $r_2$ in $(1/\sqrt{2b}, \infty)$. All these further imply that $F$ is increasing in $(0, r_1)$ and $(r_2, \infty)$ and decreasing in $(r_1, r_2)$.
Notice that $F^\prime(\sqrt{\log(2ab)/b}) = 1 - \sqrt{\log(2ab)/b} > 0$, since we have assumed that $\log(2ab) < b$, so we have that $r_2 \in (1/\sqrt{2b}, \sqrt{\log(2ab)/b})$.
Also notice that $F(1/\sqrt{2b}) = a\exp(-1/2) + 1/\sqrt{2b} < a = F(0)$, since we have assumed that $a\sqrt{b} > 1/(1 - \exp(-1/2))$. Therefore, $0$ is not the minimum point.
In summary, the minimum point of $F$ is in $(1/\sqrt{2b}, \sqrt{\log(2ab)/b})$ and satisfies $2abx\exp(-bx^2) = 1$.
\end{proof}

\section*{Acknowledgements}
The first author's research was supported in part by the NSF Grants: DMS-1561141 and DMS-1613016. The authors are sincerely grateful to the Associate Editor and two anonymous referees for their many corrections and suggestions that help to significantly improve the manuscript.}

%
%
%

\end{document}